\documentclass[mnsc]{informs3} 
\RequirePackage{bm}
\RequirePackage{endnotes}

\OneAndAHalfSpacedXII 

\usepackage{algorithm}
\usepackage{algpseudocode}
\usepackage{tikz}
\usepackage{algorithm}      
\usepackage{algpseudocode}
\usepackage{amsmath,mathtools}
\usepackage{bbm,bm}
\usepackage{siunitx}
\usepackage{hyperref}
\usepackage{tikz}
\usetikzlibrary{positioning}
\usepackage{pgfplots}
\pgfplotsset{width=5cm, compat=1.16}
\usepackage{subcaption} 
\pgfplotsset{ every non boxed x axis/.append style={x axis line style=-},every non boxed y axis/.append style={y axis line style=-}}

\newcommand{\R}{\mathbb{R}}
\newcommand{\Z}{\mathbb{Z}}
\newcommand{\E}{\mathbb{E}}
\newcommand{\data}{\mathcal{D}}

\newcommand{\hth}{\hat{\theta}}
\newcommand{\hsig}{\hat{\sigma}}
\newcommand{\Htheta}{H_\theta(\bm \theta)}

\newcommand{\hbth}{\hat{\bm{\theta}}}
\newcommand{\blam}{\bm{\lambda}}
\newcommand{\bth}{\bm{\theta}}

\DeclareMathOperator{\I}{\mathbb{I}}
\newcommand{\one}{\bm{1}} 
\DeclareMathOperator{\diag}{diag}

\usepackage{natbib}
 \bibpunct[, ]{(}{)}{,}{a}{}{,}%
 \def\bibfont{\small}%
\EquationsNumberedThrough    

\TheoremsNumberedThrough     
\ECRepeatTheorems  %

\MANUSCRIPTNO{MNSC-0001-2024.00}

\begin{document}



\RUNTITLE{Post-estimation Adjustments in Data-driven Decision-making with Applications in Pricing}

\TITLE{Post-estimation Adjustments in Data-driven Decision-making with Applications in Pricing}

\ARTICLEAUTHORS{%
\AUTHOR{Michael Albert}
\AFF{ Darden School of Business, University of Virginia, \EMAIL{albertm@darden.virginia.edu}}
\AUTHOR{Max Biggs}
\AFF{Darden School of Business, University of Virginia,\EMAIL{biggsm@darden.virginia.edu}}
\AUTHOR{Ningyuan Chen}
\AFF{Department of Management, University of Toronto Mississauga,
Rotman School of Management, University of Toronto, \EMAIL{ningyuan.chen@utoronto.ca}}
\AUTHOR{Guan Wang}
\AFF{Rotman School of Management, University of Toronto, \EMAIL{wguan.wang@rotman.utoronto.ca}}
 } 
\ABSTRACT{%
The predict-then-optimize (PTO) framework is a standard approach in data-driven decision-making, where a decision-maker first estimates an unknown parameter from historical data and then uses this estimate to solve an optimization problem. While widely used for its simplicity and modularity, PTO can lead to suboptimal decisions because the estimation step does not account for the structure of the downstream optimization problem. We study a class of problems where the objective function, evaluated at the PTO decision, is asymmetric with respect to estimation errors. This asymmetry causes the expected outcome to be systematically degraded by noise in the parameter estimate, as the penalty for underestimation differs from that of overestimation. To address this, we develop a data-driven post-estimation adjustment that improves decision quality while preserving the practicality and modularity of PTO. We show that when the objective function satisfies a particular curvature condition, based on the ratio of its third and second derivatives, the adjustment simplifies to a closed-form expression. This condition holds for a broad range of pricing problems, including those with linear, log-linear, and power-law demand models. Under this condition, we establish theoretical guarantees that our adjustment uniformly and asymptotically outperforms standard PTO, and we precisely characterize the resulting improvement. Additionally, we extend our framework to multi-parameter optimization and settings with biased estimators. Numerical experiments demonstrate that our method consistently improves revenue, particularly in small-sample regimes where estimation uncertainty is most pronounced. This makes our approach especially well-suited for pricing new products or in settings with limited historical price variation.
}%





\maketitle


\section{Introduction}\label{sec:Intro}

In data-driven decision-making, a decision-maker (DM) often relies on historical data to learn about an operational environment and make optimal choices. This process typically involves estimating unknown parameters that characterize the environment. A standard method for such problems is the predict-then-optimize (PTO) framework. In this two-stage approach, parameters are first estimated from data using statistical or machine learning models; these estimates are then used as inputs to an optimization problem, which is solved to yield a final decision.

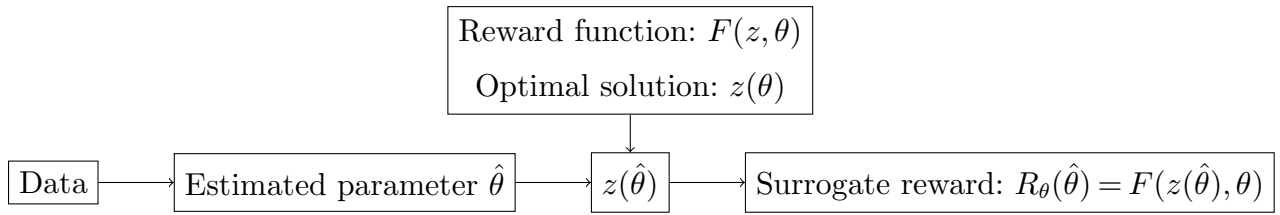
\begin{figure}[ht]
    \centering
    \begin{tikzpicture}
        \node (data) [draw, rectangle] {Data};
        \node (theta) [draw, rectangle, right=1cm of data] {Estimated parameter \(\hat{\theta}\)};
        \node (z) [right=1cm of theta, draw, rectangle] {\(z(\hat{\theta})\)};
        \node (obj) [right=1cm of z, draw, rectangle] {Surrogate reward: \(R_{\theta}(\hat \theta)=F(z(\hat \theta), \theta)\)};
        \node (opt) [above=0.5cm of z, draw, rectangle, align=center] {Reward function: $F(z, \theta)$\\ Optimal solution: $z(\theta)$};
        
        \draw[->] (data) -- (theta);
        \draw[->] (theta) -- (z);
        \draw[->] (z) -- (obj);
        \draw[->] (opt) -- (z);

    \end{tikzpicture}
    \caption{The workflow of predict-then-optimize (PTO) in data-driven decision-making}
    \label{fig:PTO_flow}
\end{figure}

To formalize this, consider a DM who seeks to choose a decision $z$ to maximize a reward function $F(z,\theta)$, where $\theta$ is an unknown parameter. We refer to $F(z,\theta)$ as the \emph{reward} function without loss of generality. If $\theta$ were known, the optimal decision would be $z(\theta)=\argmax_z F(z,\theta)$, assuming an unconstrained problem for simplicity. In the PTO framework, the DM first uses historical data to produce an estimate, $\hth$. This is the ``predict'' step. The DM then solves the optimization problem using this estimate, yielding the decision $z(\hth)$. This is the ``optimize'' step. The resulting reward is $F(z(\hth),\theta)$. Since the functional form of $F$ and thus the mapping $z(\cdot)$ are known, we can analyze the performance of PTO by studying the \emph{surrogate reward} function, defined as $R_\theta(\hth) \coloneqq F(z(\hth),\theta)$. This workflow is illustrated in Figure~\ref{fig:PTO_flow}.

The PTO framework is popular due to its simplicity and modularity. By treating estimation and optimization as distinct steps, it allows practitioners to leverage well-established statistical and machine learning methods for the prediction task. However, this separation is also its primary weakness. Recent literature (e.g., \citet{elmachtoub2022smart,chuSolvingPriceSettingNewsvendor2023}) has shown that decision quality can be improved by incorporating the structure of the downstream optimization problem into the estimation process. This integrated paradigm, often called decision-focused learning, aims to produce estimates that are not necessarily the most accurate in a statistical sense, but are the most useful for the ultimate decision. It accounts for how estimation errors propagate through the optimization problem, thereby aligning the predictive model with the DM's final objective.

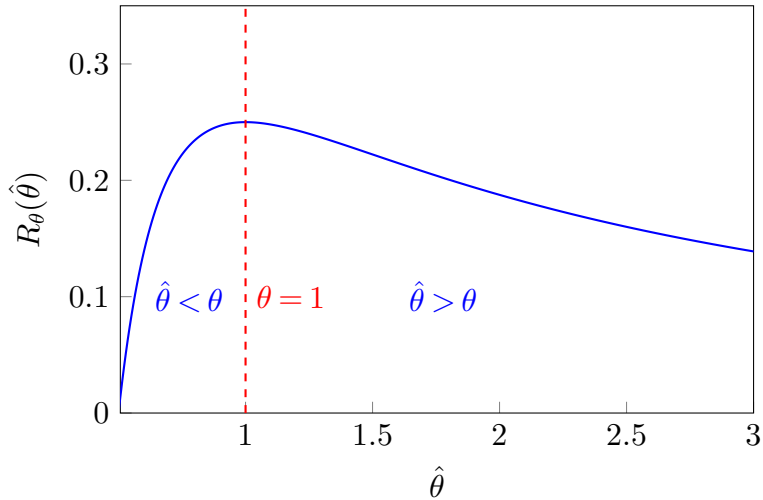
\begin{figure}[ht]
    \centering
        \begin{tikzpicture}
  \begin{axis}[
    width=10cm, height=7cm,
    xlabel = $\hat\theta$, ylabel = {$R_\theta(\hth)$},
    samples=400,
    domain=0.1:3,
    ymin=0, ymax=0.35,
    axis line style={draw=black},
    xtick pos=left,
    ytick pos=left,
    enlargelimits=false
  ]
      \addplot[color=blue,thick]{1/(2*x)*(1-1/(2*x))};
    \addplot[red, thick, dashed] coordinates {(1, 0) (1, 0.35)};
    \node[anchor=west, red] at (axis cs:1,0.1) {$\theta=1$};
    \node[anchor=west, blue] at (axis cs:0.6,0.1) {$\hat\theta<\theta$};
    \node[anchor=west, blue] at (axis cs:1.6,0.1) {$\hat\theta>\theta$};
  \end{axis}
  \end{tikzpicture}
      \caption{The surrogate reward $R_\theta(\hth)$ for the linear demand pricing problem (Example~\ref{exp:pricing}), with true parameters $a=1$ and $\theta=1$. The function is asymmetric around its maximum at $\hth=\theta$. An underestimation ($\hth < \theta$) results in a larger revenue loss than an overestimation of the same magnitude.}
    \label{fig:assymetric_revenue}
\end{figure}

Our work focuses on a specific class of problems where the PTO decision rule $z(\hth)$ has a closed-form expression. In such cases, we can analytically study how estimation errors in $\hth$ affect the final reward. We use the following data-driven pricing problem as a running example to illustrate the core concepts of our approach.
\begin{example}[Data-driven pricing with linear demand using PTO]\label{exp:pricing}
    A decision-maker (DM) sets a price $p$ to maximize revenue $F(p,\theta)=pd(p)$, where the demand function is
    \begin{equation}\label{eq:demand-func}
    d(p) = a-\theta p.
\end{equation}
For simplicity, we assume $a$ is known and the price sensitivity $\theta$ is unknown.\footnote{This basic framework can be applied to the case where an incumbent price and demand $(p_0,d_0)$ satisfying $d_0=a-\theta p_0$ are known, but $a$ and $\theta$ are not. This scenario arises when an item has been sold at a fixed price $p_0$ for an extended period, but recent price experimentation provides sufficient variation to estimate the full demand curve. In later sections, we extend our framework to multiple parameters and alternative demand functions.}
If the demand function were fully known, the optimal price would be $p(\theta) = a/(2\theta)$.
Under the PTO approach, the DM first obtains an estimate $\hth$ for the price sensitivity and then sets the price to $p(\hth) = a/(2\hth)$. The resulting surrogate reward, as a function of the estimate, is
\begin{align}\label{eq:pricing-objective}
    R_\theta(\hth) &= p(\hth) d(p(\hth)) = \frac{a}{2\hth}\left(a-\theta \frac{a}{2\hth}\right) \\
    &= \frac{a^2}{2\hth} - \frac{a^2\theta}{4\hth^2}. \notag
\end{align}
\end{example}
This example illustrates a general phenomenon: the structure of the PTO decision can introduce a critical asymmetry in the surrogate reward $R_\theta(\hth)$ with respect to the estimate $\hth$. As shown in Figure~\ref{fig:assymetric_revenue}, while the reward is maximized at the true parameter value, $\hth = \theta$, the penalty for underestimation ($\hth < \theta$) is far more severe than the penalty for an equivalent overestimation ($\hth > \theta$). Consequently, for any noisy estimator $\hth$, even an unbiased one where $\mathbb{E}[\hth] = \theta$, the expected reward $\mathbb{E}[R_\theta(\hth)]$ is reduced more by the possibility of underestimation than overestimation. This suggests that a corrective adjustment that intentionally biases $\hth$ upward could improve expected performance by mitigating the downside risk of underestimation. Notably, this adjustment stands in contrast to traditional regularization techniques like ridge regression \citep{hoerl1970ridge} or LASSO \citep{tibshirani1996regression}, which typically shrink estimates toward zero to improve predictive accuracy. Our analysis highlights the need to tailor the bias-variance trade-off to the specific structure of the downstream decision problem.

In this paper, we develop a general framework for a post-estimation adjustment that systematically improves upon the standard PTO approach. By explicitly accounting for the structure of the reward function, our method corrects the PTO estimate to improve decision quality. This correction operates as a simple, data-driven step applied after the initial prediction, preserving the modularity of the PTO workflow and allowing practitioners to enhance existing models without altering the underlying estimation procedure. Our main contributions are as follows:
\begin{itemize}
    \item We first derive an oracle adjustment that establishes how PTO can be systematically improved when the true value of the unknown parameter is known. 
    It is similar to the LASSO analysis where the hyperparameter for the $\ell_1$ penalty depends on the true coefficients \citep{tibshirani1996regression, van2009conditions}.
    This oracle result serves as a theoretical benchmark, demonstrating that adjustments to an unbiased estimated parameter lead to improved decision outcomes compared to standard PTO.  
    Additionally, we provide formal performance bounds that characterize the precise asymptotic improvement achievable through this approach. We show that this improvement is optimal relative to the class of second-order adjustments we study.
    \item We characterize the conditions under which these improvements hold, in particular, when the objective resulting from the PTO decision satisfies a condition on the ratio of the third- to the second-order derivative under mild regularity conditions on the estimation procedure. 
    When this assumption is satisfied, the Taylor expansion of the reward simplifies significantly, and a simple closed-form adjustment can be derived, for which we are able to characterize the uniform improvement of the objective. When the derivative ratio non-trivially satisfies this condition, this implies an asymmetry of the objective with respect to the estimated parameter. Surprisingly, this condition is met in many common pricing problems, and our results extend beyond the linear pricing setting to widely used log-linear and power law demand models.
    \item Building on the oracle case, we extend our framework to a data-driven setting. We propose a plug-in-adjustment that follows the same functional form of the oracle, but is learned from observed data rather than relying on knowledge of the true parameter. We prove that our method guarantees a uniform asymptotic improvement over PTO, provided that the estimation procedure is unbiased. Again, this relies on simplification due to the condition of the ratio of the third to the second derivative of the reward function, and does not hold for arbitrary functions. Counterintuitively, we show that the data-driven improvement is actually larger than the oracle improvements for many objectives.
    \item Beyond the single-parameter setting, we extend our framework to multi-parameter optimization problems. While an oracle adjustment that improves over PTO can be derived with a generalization of the necessary assumptions, the plug-in approach no longer achieves a universal data-driven improvement. However, we propose a bootstrap procedure that asymptotically recovers the oracle adjustment, and performs well empirically. 
    \item Finally, we validate our theoretical findings through extensive numerical experiments, illustrating the effectiveness of our post-estimation adjustment across a range of pricing settings. Our results confirm that the proposed method consistently improves decision quality over standard PTO, particularly in small-sample regimes, where estimation uncertainty is most pronounced, making it especially well-suited for pricing new products or those with limited price variation. 
\end{itemize}


By providing both theoretical foundations and practical considerations, this work advances the PTO framework by introducing decision-aware corrections that enhance optimization outcomes while preserving the simplicity and modularity that make PTO widely used in practice.

\section{Related literature}


There is significant recent literature on improving predict then optimize by incorporating the operational decision within the estimation of an unknown parameter, variously known as smart predict then optimize \citep{elmachtoub2022smart}, decision-focused learning \citep{mandi2024decision}, end-to-end optimization \citep{donti2017task}, or integrated learning and optimization \citep{sadana2024survey}. A prominent approach is decision rule learning \citep{ban2019big,bertsimas2022data}, where a decision rule (a mapping from the features to a decision) is inserted directly into the DM's objective and optimized using empirical risk minimization. 
In smart-predict-then-optimize \citep{elmachtoub2022smart}, a contextual decision, optimal with respect to a given estimated parameter, is compared to true optimal given the realized parameter. 
This regret is minimized to estimate the parameter, via a convex surrogate. 
End-to-end optimization \citep{donti2017task,amos2017optnet} integrates differentiable optimization layers into machine learning models, allowing the full decision pipeline to be trained jointly. 
Interestingly, \cite{hu2022fast} show that standard PTO can achieve faster rates of convergence than end-to-end approaches.
\cite{sadana2024survey} provides a comprehensive review of this literature.

Much of the existing research on improving predict-then-optimize (PTO) methods assumes that although the unknown parameter is not observed at the time of decision-making (often relying instead on contextual information), it becomes observable afterward \citep{ban2019big,elmachtoub2022smart}. For example, in stochastic shortest-path problems, as discussed in \citet{elmachtoub2022smart}, past observations of travel times across all routes enable empirical evaluation of counterfactual outcomes, which in turn inform improved decisions. 
\cite{liu2023active} propose an active learning framework in a similar class of problems.
However, in many practical applications, the parameter of interest is not directly observed but is embedded within a structural equation relating it to observable quantities. 
This challenge occurs in pricing applications, where firms seek to estimate demand functions in order to optimize prices. In such cases, key parameters governing consumer behavior, such as demand elasticity or willingness to pay, are latent, and firms only observe realized sales at given prices. This lack of direct observation complicates counterfactual evaluation and precludes the straightforward application of standard smart PTO methods.

The operational data analytics (ODA) framework \citep{feng2023framework} investigates a similar class of data-driven decision-making problems that satisfy a specific homogeneity condition, which characterizes how the decision, the unknown parameters, and the objective scale relative to each other. 
When the condition holds, ODA can provide a data-driven decision that is uniformly optimal for all the possible values of the unknown parameter.
The structure of our framework and the assumptions are conceptually related to the ODA framework. 
However, we believe the running example for linear pricing, and thus our general framework, does not fit into ODA.
In particular, in Example~\ref{exp:pricing}, if we treat the residual as the unknown distribution $X$ in \cite{feng2023framework}, then we have weaker assumptions on the problem structure.
For example, we allow the residual to be arbitrarily distributed and the focal unknown parameter to be unrelated to the distribution.
In contrast, we do not target at the uniform optimality of our approach but show the improvement relative to PTO asymptotically.
The ODA framework has been applied to various problems that satisfy this condition, such as the newsvendor problem \citep{chu2023solving, siegel2021profit}, assortment optimization \citep{feng2022consumer}, and service speed design \citep{feng2024operational}.

There is also a growing interest in data-driven pricing. Incorporating the downstream pricing decision, \citet{banPersonalizedDynamicPricing2021} and \citet{javanmard2019dynamic} study how to regularize the contextual (and potentially high-dimensional) intercept term in a pricing model to achieve optimal regret rates. In contrast, our approach focuses on the challenge of adjusting the price sensitivity term. \citet{chernozhukov2019semi} present a doubly robust estimation approach for semi-parametric models that uses contextual pricing as an example. In this case, the robustness includes being robust to misspecification of the estimated demand function. Similar to our approach, which also aims to price effectively with few samples, \citet{fu2015randomization,babaioff2018two,daskalakis2020more,allouah2022pricing}, and  \citet{allouah2023optimal} establish worst-case revenue bounds based on limited price-demand observations. Further examples of robust pricing include \citet{guan2023randomized} and \citet{chen2022distribution}. In comparison, we look at average case performance, which can be less conservative and more practically relevant in many settings. \citet{ito2018unbiased} highlights how the optimizer's curse \citep{smith2006optimizer} can occur in pricing, whereby an optimized objective, based on an observed uncertain parameter, systematically overestimates the true optimal objective. While \citet{ito2018unbiased} provide a methodology for adjusting for this, by analyzing specific forms of demand, we can exactly characterize an adjustment in closed form and precisely characterize the rate of improvement over PTO. We refer to \citet{chen2023frontiers} for a comprehensive recent review of data-driven pricing.  

Our work is also broadly related to the statistics literature on regularization or shrinkage of parameter estimates \citep{tibshirani1996regression,hoerl1970ridge}. Also related is the literature on debiasing non-linear estimates, for example, Stein's Unbiased Risk Estimate (SURE) \citep{stein1981estimation}. 
\cite{loke2022decision} propose to combine regularization and PTO for the data-driven decision-making problem studied in \cite{elmachtoub2022smart}, which does not cover pricing.
We extend these ideas by demonstrating how debiasing techniques can be tailored to a specific class of functions, in particular, those that emerge in PTO-based pricing settings, and rigorously characterizing the resulting performance improvements.

\section{Framework for Single-parameter PTO Adjustment}\label{sec:model}
We begin by analyzing the surrogate reward function $R_{\theta}(\hat{\theta}) = F(z(\hat{\theta}),\theta)$, which is obtained by substituting the PTO decision $z(\hat{\theta})$ into the true objective $F(z,\theta)$ (see Figure~\ref{fig:PTO_flow}). While the functional form of $F$ is known, it depends on an unknown parameter $\theta$. The PTO estimator, $\hat{\theta}$, is typically produced by a standard statistical or machine learning model that is agnostic to the structure of the downstream reward function $R_{\theta}(\cdot)$.

Our central idea is to apply a \emph{post-estimation adjustment} to the PTO estimator, transforming $\hat{\theta}$ into a new estimate $\pi(\hth)$ that yields a higher expected reward, i.e., $\E[R_{\theta}(\pi(\hth))] > \E[R_{\theta}(\hth)]$, where the expectation is over the sampling distribution of $\hth$. Since the true parameter $\theta$ is unknown, this adjustment function $\pi(\cdot)$ cannot depend on it. This raises the central question of our work: Is it possible to construct an adjustment that improves the expected reward \emph{uniformly} across all possible values of $\theta$? This section establishes a framework in which the answer is yes.

\paragraph{Setup and Notation.}
We consider a sequence of estimators $\hat{\theta}_n$ generated from datasets of increasing size $n$. We impose the following mild regularity conditions, which are satisfied by a wide class of common estimators.
\begin{assumption}[Regularity of the Estimator]\label{asp:eps-dist}
As $n\to\infty$, we assume:
\begin{enumerate}
    \item The estimator is root-$n$ consistent: The mean squared error $\sigma_n^2 \coloneqq \E[(\hth_n-\theta)^2]$ satisfies $\lim_{n\to \infty} n\sigma_n^2 = \sigma_\theta^2$ for some constant $\sigma_\theta^2 > 0$.
    \item The estimator has a diminishing bias: $\E[|\hth_n-\theta|]=o(n^{-1})$.
    \item Higher-order centered moments are well-behaved: $\E[|\hth_n-\theta|^k]=O(n^{-k/2})$ for $k=3,4,\dots$.
\end{enumerate}
\end{assumption}
Assumption~\ref{asp:eps-dist}(1) restricts our focus to estimators that are nearly \emph{unbiased}; if a first-order bias were to dominate the sampling error, a direct debiasing step would be more effective than our second-order adjustment. Parts (2) and (3) are standard conditions for estimators that converge at a rate of $1/\sqrt{n}$, ensuring that their moments scale appropriately. To make this concrete, we consider the OLS estimator for our running example.
\paragraph{Estimation of $\hth$ in Example~\ref{exp:pricing} using OLS.}
Given a dataset of price and demand pairs $\mathcal{D}\coloneqq\{p_i,d_i\}_{i=1}^n$ from the model
\begin{equation}\label{eq:data-generation}
  d_i= a - \theta p_i +\epsilon_i,
\end{equation}
where the errors $\epsilon_i$ are i.i.d.\ with $\E[\epsilon_i]=0$ and finite higher-order moments, the ordinary least squares (OLS) estimator for $\theta$ is
\begin{equation}\label{eq:ols-pricing}
    \hth= -\frac{\sum_{i=1}^n p_i(d_i-a)}{\sum_{i=1}^n p_i^2}.
\end{equation}
Under standard conditions on the covariates $p_i$, this estimator can be shown to satisfy Assumption~\ref{asp:eps-dist}.

Next, we formalize our assumptions on the surrogate reward function $R_{\theta}(\cdot)$. Let $R^{(k)}_{\theta}(x)$ denote the $k$-th derivative of $R_{\theta}(x)$ with respect to its argument $x$.

\begin{assumption}[Structure of the Surrogate Reward]\label{asp:r-derivative}
For any true parameter $\theta$ in its domain, we assume:
\begin{enumerate}
\item The reward is maximized at the true parameter value: $R^{(1)}_{\theta}(\theta)=1$ and $R^{(2)}_{\theta}(\theta)<0$.
\item The ratio of the third to the second derivative is proportional to $1/\theta$: there exists a constant $C$, independent of $\theta$, such that $\frac{R^{(3)}_{\theta}(\theta)}{R^{(2)}_{\theta}(\theta)}= \frac{C}{\theta}$.
\item The fifth derivative is locally bounded: $|R^{(5)}_{\theta}(x)| \le M$ for some constant $M$ and for all $x$ on the interval between $\theta$ and $\hat{\theta}_n$.
\end{enumerate}
\end{assumption}
Part (1) is a natural condition stating that the reward is locally maximized when the parameter estimate is perfectly accurate. 
Part (2) is the crucial structural assumption of our framework. While it may seem restrictive, it is satisfied by the surrogate reward functions derived from many common economic models, as verified in Proposition~\ref{prop:function-2-3derivative}. 
Part (3) is a mild technical requirement for bounding the error terms in our Taylor series analysis. 
This condition implicitly assumes the existence of all lower-order derivatives up to order five.
It bounds the fifth derivative only along the data-dependent segment between $\theta$ and $\hat{\theta}_n$, a much weaker requirement than a global bound and easily satisfied after practical truncation (e.g. enforcing $\hat{\theta}_n>\epsilon$ in Example~\ref{exp:pricing}).
The next result demonstrates the generality of the assumption.
\begin{proposition}\label{prop:function-2-3derivative}
    Assumption~\ref{asp:r-derivative}(2) is satisfied for the following classes of surrogate reward functions $R_\theta(\hth)$:
    \begin{itemize}
        \item $R_\theta(\hth)=f(\theta)\sum_{i=1}^I a_i \theta^{j_i}\hth^{K-j_i}+g(\theta) \hth+h(\theta)$ for $I,K\in \Z_+$, $j_i\in \R$.
        \item $R_\theta(\hth)=f(\theta)\exp\left( \frac{a\hth}{\theta}\right)$ or $R_\theta(\hth)=\frac{f(\theta)}{\hth}\exp\left( \frac{a\theta}{\hth}\right)$ for $a\in \R$.
        \item $R_\theta(\hth)=f(\theta)\log\left( g(\theta)\hth\right)+\hth$ for $g(\theta)\hth >0$.
    \end{itemize}
   Furthermore, it is satisfied for the following objective functions $F(z, \theta)$, when $z(\hth)$ is derived from the first-order condition:
    \begin{itemize}
     \item $F(z, \theta) = a_1 \theta \log(z)-a_2z$, for $a_1\theta>0, a_2>0${}
     \item $F(z, \theta) = f(a,\theta, \gamma)z(a-\theta z)^{\gamma}$, for $\theta>0$ and known $\gamma>0, a>0$.
     \item $F(z, \theta) = f(a,\theta)z\exp(a-\theta z)$, for $a>0, \theta>0$.
     \end{itemize}
\end{proposition}
To demonstrate the practical relevance of Assumption~\ref{asp:r-derivative}, we verify it for our running example and two other widely used pricing models. For the linear demand model in Example~\ref{exp:pricing}, the surrogate reward is $R_\theta(\hth)=-a^2\theta/(4\hth^2)+a^2/(2\hth)$, for which it is straightforward to verify that $C=-6$. We now consider two other canonical demand models.
\begin{example}[Pricing with Log-Linear Demand]\label{exp:log-linear}
Consider a firm setting a price $p$ to maximize revenue where demand is log-linear \citep{talluri2006theory}: $d(p) = \exp(a-\theta p)$. Assuming $a$ is known, the revenue is $F(p,\theta)= p\exp(a-\theta p)$, and the optimal price is $p(\theta)=1/\theta$. The PTO price is $p(\hat{\theta})=1/\hat{\theta}$, which yields the surrogate reward
\begin{equation}
\label{eq:log-linear-pricing-objective}
    R_\theta(\hth) = p(\hth) d(p(\hth)) = \frac{1}{\hth}\exp\left(a- \frac{\theta}{\hth}\right).
\end{equation}
This functional form is covered by Proposition~\ref{prop:function-2-3derivative}. By direct differentiation, we find
\begin{equation}
\label{eq:log-linear-c}
\frac{R_\theta^{(3)}(\theta)}{R_\theta^{(2)}(\theta)}=-\frac{4}{\theta},
\end{equation}
and thus Assumption~\ref{asp:r-derivative}(2) holds with $C=-4$.
\end{example}

\begin{example}[Pricing with Power-Law Demand]\label{exp:power-law}
Now suppose demand follows a power law \citep{caplin1991aggregation}, $d(p)= (a-\theta p)^\gamma$, where $\gamma>0$ is a known parameter. The revenue is $F(p,\theta)= p(a-\theta p)^\gamma$, the optimal price is $p(\theta)= a/(\theta(1+\gamma))$, and the surrogate reward is
\[R_\theta(\hat{\theta})
=\frac{a^{\gamma+1}}{(1+\gamma)\hat{\theta}}\left(1-\frac{\theta}{\hat{\theta}(1+\gamma)}\right)^\gamma.
\]
Differentiating with respect to $\hth$ and evaluating at $\hth=\theta$ yields
\[
R^{(2)}_\theta(\theta)=-\frac{a^{\gamma+1}}{\gamma}\left(\frac{\gamma}{\gamma+1}\right)^{\gamma}\frac{1}{\theta^3} \quad \text{and} \quad R^{(3)}_\theta(\theta)=-\frac{a^{\gamma+1}}{\gamma^2}\left(\frac{\gamma}{\gamma+1}\right)^{\gamma}\frac{2(1+2\gamma)}{\theta^4}.
\]
Therefore, the ratio is given by
\[
\frac{R_\theta^{(3)}(\theta)}{R_\theta^{(2)}(\theta)}
= -\frac{2(1+2\gamma)}{\gamma}\frac{1}{\theta},
\]
which satisfies Assumption~\ref{asp:r-derivative}(2) with $C = -2(1+2\gamma)/\gamma$.
\end{example}

It is instructive to note that standard statistical loss functions also fit within this framework. For example, the negative mean squared error, $R_\theta(\hth)=-(\theta-\hth)^2$, satisfies Assumption~\ref{asp:r-derivative} with $C=0$. As we will show, our adjustment in this case reduces to a form of shrinkage similar to ridge regression. This connection highlights that our framework provides a decision-aware generalization of the classical bias-variance trade-off. The following sections use these assumptions to develop an adjustment that systematically improves upon PTO.


 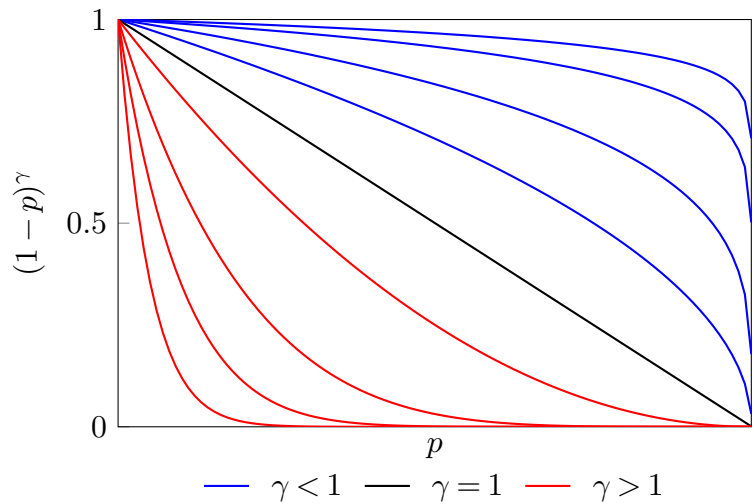
\begin{figure}[htbp] 
 \centering
    \begin{tikzpicture}
  \begin{axis}[
    width=10cm, height=7cm,
    xlabel = $p$, ylabel = {$ (1 - p)^\gamma$},
    domain=0:0.999,
    samples=100,
    ymin=0, ymax=1,
    axis line style={draw=black},
    xtick pos=left,
    ytick pos=left,
    enlargelimits=false,
    xtick=\empty,
    ytick={0,0.5,1}
  ]
      \addplot[color=blue,thick] {(1 - x)^0.05};\label{plot:g1}
      \addplot[ color=blue,thick]{(1 - x)^0.1};
      \addplot[color=blue,thick] {(1 - x)^0.25};
      \addplot[ color=blue,thick]{(1 - x)^0.5};
      \addplot[color=black,thick]{(1 - x)};\label{plot:g2}
      \addplot[ color=red,thick] {(1 - x)^2};\label{plot:g3}
      \addplot[color=red,thick]  {(1 - x)^5};
      \addplot[color=red,thick]  {(1 - x)^10};
      \addplot[color=red,thick]  {(1 - x)^20};
  \end{axis}
\end{tikzpicture}

  \begin{tikzpicture}[overlay, remember picture]
    \node[anchor=south, xshift=0.75cm, yshift=0cm] {
      \begin{tabular}{llllll}
        \ref{plot:g1} & $\gamma<1$ &
        \ref{plot:g2} & $\gamma=1$ &
        \ref{plot:g3} & $\gamma>1$
      \end{tabular}
    };
  \end{tikzpicture}

\caption{The flexibility of the power-law demand model from Example~\ref{exp:power-law}. The plot shows demand curves $d(p)=(1-p)^\gamma$ for various values of the exponent $\gamma$. For $\gamma=1$ (solid black line), demand is linear. For $\gamma > 1$ (red lines), demand is convex. For $\gamma < 1$ (blue lines), demand is concave.}
\end{figure}



\subsection{The Oracle Adjustment}\label{sec:oracle-adj}
To establish a theoretical benchmark, we first consider an idealized \emph{oracle} setting where the decision-maker knows both the true parameter $\theta$ and the asymptotic variance of its estimator, $\sigma_\theta^{2} = \lim_{n\to\infty} n\,\mathrm{Var}(\hat{\theta}_n)$. We analyze the effect of a simple multiplicative adjustment of the form
\[
\pi(\hth_n)=\hth_n\left(1+\frac{\lambda}{n}\right),
\]
where $\lambda$ is an adjustment coefficient that can depend on the oracle information ($\theta$ and $\sigma_\theta$) but not on the sample size $n$. Setting $\lambda=0$ recovers the standard PTO decision. The oracle's problem is to choose the $\lambda$ that maximizes the expected reward:
\[
   \max_{\lambda\in\mathbb{R}}
   \E\left[
        R_\theta\left(\hth_n\left(1+\frac{\lambda}{n}\right)\right)
      \right].
\]
The following proposition provides a closed-form solution for the optimal $\lambda$ and characterizes the resulting improvement in performance.
\begin{proposition}[Oracle Adjustment]\label{prop:known-beta}
Suppose Assumptions~\ref{asp:eps-dist} and \ref{asp:r-derivative} hold. Let the oracle adjustment coefficient be
\begin{equation}\label{eq:single-par-opt-lambda}
\lambda^\ast = -\frac{(C+2)\sigma^2_\theta}{2\theta^2},
\end{equation}
where $C$ is the constant from Assumption~\ref{asp:r-derivative}(2). Then the asymptotic improvement in expected reward over PTO is given by
\begin{equation}\label{eq:single-par-improvement-gap}
\lim_{n\to\infty} n^2\left(\E\left[R_{\theta}\left(\hth_n\left(1+\frac{\lambda^\ast}{n}\right)\right)\right]-\E[R_{\theta}(\hth_n)]\right) = -\frac{R^{(2)}_{\theta}(\theta)(C+2)^2\sigma_\theta^4}{8\theta^2}\ge 0.
\end{equation}
Moreover, this choice of $\lambda^\ast$ is optimal in the sense that the coefficient $-(C+2)/2$ is the unique constant multiplying $\sigma_\theta^{2}/\theta^{2}$ that maximizes the asymptotic improvement for all valid $(\theta,\sigma_\theta)$.
\end{proposition}

We now discuss the implications of this result.
\paragraph{The Form of the Adjustment.}
The multiplicative factor $(1+\lambda/n)$ is chosen to isolate the second-order effects of estimation error. A Taylor expansion of the expected reward reveals that the leading-order terms affected by $\lambda$ are of magnitude $O(n^{-2})$ and form a simple quadratic function of $\lambda$. This structure allows for a straightforward optimization to find the optimal $\lambda^\ast$. Assumption~\ref{asp:r-derivative}(2) is critical here; it ensures that the optimal coefficient depends only on the known constant $C$, rather than on the specific values of the derivatives $R_\theta^{(2)}(\theta)$ and $R_\theta^{(3)}(\theta)$. This simplification is the key to developing a practical, data-driven version of the adjustment in the next section.

\paragraph{Magnitude and Drivers of the Improvement.}
The improvement from the adjustment is of order $O(n^{-2})$. While this may seem small, it is the appropriate scale for a second-order correction applied to a nearly unbiased, root-$n$ consistent estimator. 
In a closely related work, \cite{gotoh2025data} demonstrate that distributionally robust optimization can outperform the sample average approximation, which is conceptually related to PTO, by the same magnitude $O(n^{-2})$ asymptotically.
Their analysis also relies on local perturbations.
However, their formulation and ours do not cover each other. 
The practical impact of this improvement can be significant, especially in small-sample settings. This is common when pricing new products \citep{baardman2018leveraging, cohen2021simple}, where historical data is limited, or when segmenting customers into many distinct groups, leaving few observations per group. The magnitude of the gain, as shown in \eqref{eq:single-par-improvement-gap}, is driven by two main factors: the variance of the estimator ($\sigma_\theta^4$) and the curvature of the reward function ($-R^{(2)}_{\theta}(\theta)$). Problems with noisier estimators or more sharply peaked reward functions stand to benefit the most from the adjustment.

\paragraph{The Direction of the Adjustment.}
The constant $C$ determines the character of the optimal adjustment. If $C=-2$, then $\lambda^\ast=0$, and the standard PTO approach is already asymptotically optimal. The sign of $C+2$ dictates the direction of the adjustment:
\begin{itemize}
    \item \textbf{Shrinkage ($C > -2$):} The adjustment shrinks the estimate toward zero ($\lambda^\ast < 0$). This aligns with conventional regularization. For instance, in a standard prediction task with negative mean squared error loss, $R_\theta(\hth)=-(\theta-\hth)^2$, we have $C=0$. Our adjustment then reduces to a form of ridge regression \citep{vanwieringenLectureNotesRidge2021}.
    \item \textbf{Expansion ($C < -2$):} The adjustment expands the estimate away from zero ($\lambda^\ast > 0$). This is a less common feature, typically unique to decision-focused objectives. As seen in our pricing examples (e.g., $C=-6$ for linear demand), applying standard shrinkage would be counterproductive and strictly worsen performance.
\end{itemize}
This analysis highlights the risks of naively applying prediction-focused tools like ridge regression to decision problems. The optimal adjustment must be tailored to the specific structure of the downstream task.

The oracle analysis provides a clear benchmark: an adjustment with the coefficient from \eqref{eq:single-par-opt-lambda} extracts all available second-order improvements. The challenge now is to replicate this gain in a practical, data-driven manner without knowledge of $\theta$ or $\sigma_\theta$.

\subsection{Data-Driven Plug-in Adjustment}\label{sec:dd-adj}
The oracle adjustment in Proposition~\ref{prop:known-beta} provides a theoretical benchmark but is impractical, as it requires knowledge of the true parameters $\theta$ and $\sigma_\theta^2$. We now develop a fully data-driven procedure that achieves the same order of improvement by replacing these unknown quantities with their sample estimates.

This ``plug-in'' approach motivates an adjustment coefficient of the form
\begin{equation}\label{eq:lambda-plugin}
    \lambda_n = k \frac{\hat{\sigma}_\theta^2}{\hat{\theta}_n^2},
\end{equation}
where $k$ is a constant to be optimized, $\hat{\theta}_n$ is the original PTO estimator, and $\hat{\sigma}_\theta^2$ is a consistent estimator for the asymptotic variance parameter $\sigma_\theta^2$. Such a variance estimate is a standard output in most statistical software packages (e.g., the squared standard error of a regression coefficient). For this plug-in approach to be well-behaved, we require the following standard assumptions.
\begin{assumption}[Regularity of the Variance Estimator]\label{asp:sigma-estimator}
The estimator $\hat{\sigma}_\theta^2$ for the asymptotic variance parameter $\sigma_\theta^2$ satisfies:
\begin{enumerate}
    \item It is nearly unbiased: $\E[\hat{\sigma}_\theta^2 - \sigma_{\theta}^2] = o(n^{-1})$.
    \item Its higher moments are well-behaved: $\E[(\hat{\sigma}_\theta^{2}-\sigma_\theta^{2})^{j}] = O(n^{-j})$ for $j=2,3,\ldots$.
    \item It is asymptotically uncorrelated with the primary estimation error: $\E[(\hth_n-\theta)(\hat{\sigma}_\theta^2-\sigma_{\theta}^2)]=o(n^{-1})$.
\end{enumerate}
\end{assumption}
\begin{assumption}[Estimator Bounded Away from Zero]\label{asp:betahat-away-from-zero}
There exists a constant $\varepsilon>0$ such that $\Pr\{|\hat{\theta}_n| \ge \varepsilon\}=1$ for all~$n$.
\end{assumption}
Assumption~\ref{asp:sigma-estimator} holds for many standard variance estimators under typical regularity conditions, i.e. linear regression. Assumption~\ref{asp:betahat-away-from-zero} prevents the denominator in~\eqref{eq:lambda-plugin} from vanishing and is easily enforced in practice by truncating $|\hat{\theta}_n|$ at a small positive threshold.
With this setup, we can now state the main result for the data-driven adjustment, which identifies the optimal constant $k$ and characterizes the resulting performance gain.

\begin{proposition}[Data-driven Plug-in Adjustment]\label{prop:unknown-beta}
Suppose Assumptions~\ref{asp:eps-dist} through \ref{asp:betahat-away-from-zero} hold. If we choose the adjustment coefficient
\begin{align*}
\lambda_n^\ast = \frac{(2-C)\hat{\sigma}^2_\theta}{2\hat{\theta}_n^2},
\end{align*}
then the asymptotic improvement in expected reward over PTO is given by
\begin{equation}\label{eq:single-par-dd-improvement}
\lim_{n\to\infty} n^2 \left(\E\left[R_{\theta}\left(\hth_n\left(1+\frac{\lambda_n^\ast}{n}\right)\right)\right]-\E[R_{\theta}(\hth_n)]\right) = -\frac{R^{(2)}_{\theta}(\theta)(2-C)^2\sigma_\theta^4}{8\theta^2}\ge 0.
\end{equation}
\end{proposition}
The proof (see Appendix~\ref{app:single-par-proof}) again relies on a Taylor expansion of the expected reward. The analysis is more intricate than in the oracle case because the adjustment coefficient $\lambda_n^\ast$ is now a random variable. Expanding this term introduces additional complexities. However, the specific functional form of our adjustment in \eqref{eq:lambda-plugin}, combined with the structural property of the reward function in Assumption~\ref{asp:r-derivative}(2), causes these complex terms to simplify, making it possible to solve for the optimal coefficient $k=(2-C)/2$.

\paragraph{Magnitude of Improvement: The Data-Driven Advantage.}
The data-driven adjustment achieves the same $O(n^{-2})$ order of improvement as the oracle. However, comparing the improvement terms in \eqref{eq:single-par-dd-improvement} and \eqref{eq:single-par-improvement-gap} reveals a surprising result: the constant $(C+2)^2$ from the oracle case is replaced by $(2-C)^2$. For many problems, including our pricing examples where $C<0$ (e.g., $C=-6$ for linear demand), we have $(2-C)^2 > (C+2)^2$; for instance, the data-driven improvement is four times the oracle improvement in Example~\ref{exp:pricing}. This implies that the data-driven adjustment yields a \emph{larger} asymptotic improvement than the oracle. This phenomenon, where using estimated components can improve performance over an oracle that uses true components in a suboptimal rule, mirrors findings in other areas of statistics, such as the efficiency gains from using estimated weights in inverse propensity weighting schemes \citep{hirano2003efficient}. The randomness in the plug-in adjustment beneficially counteracts the estimation error in the PTO decision.

\paragraph{Direction of Adjustment.}
The optimal coefficient for the data-driven adjustment, $k=(2-C)/2$, also differs from its oracle counterpart, $k=-(C+2)/2$. This change accounts for the randomness in the plug-in term $\lambda_n^\ast$ and leads to another striking result: for any $C \in (-2, 2)$, the oracle and data-driven adjustments move in opposite directions. Consider the case of a standard prediction loss (MSE), where $C=0$. The oracle adjustment shrinks the estimate ($\lambda^\ast < 0$), as expected from regularization theory. In contrast, the optimal data-driven adjustment \emph{expands} it ($\lambda_n^\ast > 0$). This counter-intuitive reversal, which to our knowledge has not been documented previously, underscores how the interplay between estimation error and decision-making can lead to unexpected optimal strategies.

In summary, we have shown that a simple, data-driven plug-in rule can systematically improve upon the standard PTO approach. We will demonstrate the practical impact of both the oracle and data-driven adjustments in Section~\ref{sec:case_study_linear_pricing} before extending this framework to the more general multi-parameter setting.

\section{Multi-parameter Adjustment}\label{sec:multi-para-oracle}
Thus far we have focused on a single unknown parameter.  
Decision problems, however, often involve a vector $\bth=(\theta_1,\dots,\theta_m)\in\R^{m}$.  
Extending the oracle logic of Section~\ref{sec:oracle-adj} to $m>1$ raises two conceptual hurdles:
\begin{enumerate}
\item the adjustment must now be a \emph{matrix} $\Lambda$ rather than a scalar, and
\item higher-order derivatives of $R_{\bm\theta}(\cdot)$ are tensors rather than matrices or vectors.
\end{enumerate}
To keep notation compact we write $x^{\otimes k}\in \R^{m\times m\times \dots\times m}$ for the $k$-fold outer (tensor) product of a vector $x\in\R^{m}$. 
In particular, $A_{i_1i_2\dots i_k} = x_{i_1}x_{i_2}\dots x_{i_k}$.
Similarly, let $A\otimes B$ be the Kronecker product of two vectors and thus $(x\otimes y)_{ij}=x_i\,y_j$.

Let $\hbth_n$ be the PTO estimator computed from a sample of size~$n$.
We make assumptions on the estimator in parallel to Assumption~\ref{asp:eps-dist}.
\begin{assumption}[Regular Estimator $m>1$]
\label{asp:multi-eps-dist}
For every $\bm\theta\in\R^{m}$:
\begin{enumerate}
\item {Diminishing bias:}
      $\|\E[\hat{\bm\theta}_n-\bm\theta]\|_\infty=o(n^{-1})$. Here $\|x\|_{\infty}=\max_i |x_i|$ for $x\in \R^m$.
\item {Root-$n$ covariance:}
      $\Sigma_n:=\E\bigl[(\hat{\bm\theta}_n-\bm\theta)
                         (\hat{\bm\theta}_n-\bm\theta)^{\!\top}\bigr]$
      satisfies $\lim_{n\to\infty}n\Sigma_n=\Sigma_\theta$.
\item {Higher moments:}
      $\E\left[\|(\hat{\bm\theta}_n-\bm\theta)^{\otimes k}\|_{\max}\right]
      =O(n^{-k/2})$ for all integers $k\ge3$. 
      Here the tensor norm $\|A\|_{\max}\coloneqq \max_{i_1,i_2,\dots}|A_{i_1i_2\dots}|$ is the maximal absolute value of the tensor elements. 
\end{enumerate}
\end{assumption}
We denote by $H_{\theta}(\bth)=\nabla^{2}_{\phi}R_{\bm\theta}(\phi)\rvert_{\phi=\bm\theta}$ the Hessian matrix of the surrogate reward function at the true parameter and by $\nabla^{k}R_{\bm\theta}(\bm\theta)$ the $k$-th derivative tensor.
We then impose the assumptions for the objective function.
\begin{assumption}[Surrogate Reward Function $m>1$]\label{asp:r-multi-derivative}
  We assume that for all parameters of interest $\bm \theta\in \R^m$:
  \begin{enumerate}
      \item The gradient of $R_{\bm\theta}(\cdot)$ is $\bm 0$ at $\bm\theta$.
      \item There exists a matrix $M\in \R^{m\times m}$ such that for any diagonal matrix $\Lambda \in \R^{m\times m}$ and $\bm\theta$, we have
      \begin{align*}   
      \left\langle  \Htheta, \Lambda M\right\rangle=&\left\langle \nabla^3 R_{\bm{\theta}}(\bm{\theta}), (\Lambda \bm \theta) \otimes \Sigma_{\theta}\right\rangle.
    \end{align*}
   \item The fifth-order derivative $\|\nabla^5R_\theta(\cdot)\|_{\max}$ is uniformly bounded.
  \end{enumerate}
\end{assumption}
Assumption~\ref{asp:r-multi-derivative}(2) is the multi-parameter analogue of the ratio condition in Assumption~\ref{asp:r-derivative}(2). It imposes a specific relationship between the second and third derivatives of the reward function, which is crucial for simplifying the optimal adjustment problem. 
In particular, when $m=1$, we can connect the two conditions by $M=C\sigma_\theta^2$.

Under these two assumptions, we adjust the estimated parameter $\hbth$ with a matrix $\Lambda \in\R^{m\times m}$.
As a result, the objective function after the adjustment is $R_{\bm \theta}((\I+\Lambda/n)\hbth_n)$, where $\I$ is an identity matrix of dimension $m$.
Our goal is to characterize the improvement compared to PTO: $\E[R_{\bm \theta}((\I+\Lambda/n)\hbth_n)]-\E[R_{\bm\theta}(\hbth_n)]$ with a properly chosen $\Lambda$.
In particular, we focus on diagonal matrices $\Lambda=\diag(\bm\lambda)$.
\begin{proposition}[Multi-parameter Oracle Adjustment]\label{prop:multi-oracle}
Suppose Assumptions~\ref{asp:multi-eps-dist} and~\ref{asp:r-multi-derivative} hold. For a diagonal adjustment matrix $\Lambda=\diag(\blam)$, the asymptotic improvement in expected reward is a quadratic function of $\blam$:
\begin{align}\label{eq:multi-improvement-gap}
\lim_{n\to\infty} n^2\left(\E[R_{\bm \theta}((\I+\Lambda/n)\bm \hth_n)]-\E[R_{\bm \theta}(\bm \hth_n)]\right) = \frac{1}{2}\left(\blam^\top A(\bth) \blam+b(\bth)^\top\blam\right),
\end{align}
where the matrix $A(\bth)\in \R^{m\times m}$ and vector $b(\bth)\in \R^m$ are given by
\begin{align}
    A(\bth)_{ij} &= (H_{\bth})_{ij}\theta_i \theta_j, \label{eq:Atheta} \\
    b(\bth) &= (H_{\bth}\odot (2\Sigma_\theta+M_\theta) )\one_m. \label{eq:btheta}
\end{align}
Here, $\odot$ denotes the element-wise Hadamard product and $\one_m$ is a vector of ones. If $A(\bth)$ is negative definite, the unique optimal adjustment is $\blam^\ast = -\frac{1}{2}A(\bth)^{-1} b(\bth)$.
\end{proposition}
Similar to the single-parameter case \eqref{eq:single-par-improvement-gap}, the improvement gap is quadratic in $\bm\lambda$. 
Moreover, the gap has an order of magnitude $O(n^{-2})$. 
Different from the single-parameter case, the coefficient matrix $A(\bth)$ can be degenerate, as demonstrated by the linear pricing example presented shortly after.
It does not invalidate the effectiveness of the adjustment; rather, there are infinite choices of $\bm\lambda$ that can attain the optimal improvement gap.
Furthermore, the quadratic form of \eqref{eq:multi-improvement-gap} doesn't depend on Assumption~\ref{asp:r-multi-derivative}(2). 
However, the form of $A$ and $b$ can be significantly simplified with the assumption.

\begin{example}[Linear Demand with Two Unknown Parameters]\label{exp:linear-demand-two-par}
We illustrate the multi-parameter framework by extending our pricing example. Suppose demand is $d(p) = \theta_1 - \theta_2 p$, where both the intercept $\theta_1$ (market size) and slope $\theta_2$ (price sensitivity) are unknown. The optimal price is $p(\bth) = \theta_1 / (2\theta_2)$, and the PTO decision is $p(\hbth) = \hat{\theta}_1 / (2\hat{\theta}_2)$. The surrogate reward is
\begin{equation*}
    R_{\bth}(\hbth) = p(\hbth) \left(\theta_1 - \theta_2 p(\hbth)\right) = \frac{\theta_1 \hat{\theta}_1}{2 \hat{\theta}_2} - \frac{\theta_2 \hat{\theta}_1^2}{4 \hat{\theta}_2^2}.
\end{equation*}
This model can be shown to satisfy Assumptions~\ref{asp:multi-eps-dist} and \ref{asp:r-multi-derivative}. To apply Proposition~\ref{prop:multi-oracle}, we compute the components of the oracle adjustment. The Hessian of $R_{\bth}(\hbth)$ evaluated at $\hbth=\bth$ is
\begin{equation}\label{eq:hessian-linear-demand}
    H_{\bth} = \frac{1}{2\theta_2} \begin{bmatrix} -1 & \theta_1/\theta_2 \\ \theta_1/\theta_2 & -\theta_1^2/\theta_2^2 \end{bmatrix}.
\end{equation}
From \eqref{eq:Atheta}, the matrix $A(\bth)$ is
\begin{equation}\label{eq:a-lnear}
    A(\bth) = \frac{\theta_1^2}{2\theta_2} \begin{bmatrix} -1 & 1 \\ 1 & -1 \end{bmatrix}.
\end{equation}
This matrix is singular, which is expected: the decision $p(\hbth)$ depends only on the ratio $\hat{\theta}_1/\hat{\theta}_2$. Adjustments to $\hat{\theta}_1$ and $\hat{\theta}_2$ that preserve this ratio (to first order) do not change the decision, leading to a family of optimal adjustments. The vector $b(\bth)$ can also be computed via \eqref{eq:btheta} (see Appendix for the derivation of $M_\theta$ and the full expression for $b(\bth)$).

To find an optimal adjustment, we solve the linear system $A(\bth)\blam = -b(\bth)/2$. Since $A(\bth)$ is singular, there are infinitely many solutions. A simple approach is to set one component to zero, say $\lambda_2=0$, and solve for the other. This yields a particular oracle solution:
\begin{equation}\label{eq:foc-linear}
    \lambda_1^\ast = - \frac{\Sigma_{\theta,11}}{\theta_1^2} + \frac{3\Sigma_{\theta,12}}{\theta_1\theta_2} - \frac{2\Sigma_{\theta,22}}{\theta_2^2}, \quad \lambda_2^\ast = 0.
\end{equation}
Plugging this solution back into \eqref{eq:multi-improvement-gap} gives the optimal improvement, demonstrating that a significant, quantifiable gain is achievable even in this more complex setting.
\end{example}

\subsection{Bootstrap for Data-driven Adjustment}\label{sec:bootstrap}
For a single parameter ($m=1$) the oracle adjustment $\lambda(\theta)$ is a scalar and Section~\ref{sec:dd-adj} shows that 
the data-driven adjustment can be adapted by simply plugging in $\lambda(\hth)$ and then finding the right coefficient.
In contrast, when $m>1$ the oracle matrix $\Lambda(\theta)$ depends on~$\theta$ through high‑order derivatives of the risk and through the sampling covariance of $\hat\theta$.  
Using a plug-in $\Lambda(\hbth)$ and then optimizing the coefficients does not yield a positive improvement gap for all $\theta$, because optimizing the coefficients inevitably leads to dependence on $\bth$.
More precisely, it is unclear whether we can conjecture a form $\Lambda(\hbth)$ and then optimize the coefficients in the improvement gap without re‑introducing unknown quantities $\bth$.
In this section, we propose a bootstrap-based approach, different from Section~\ref{sec:dd-adj}, which leads to a data-driven adjustment matrix $\hat\Lambda$.

Let $\mathcal D_n(\bth)$ denote the size‑$n$ sample drawn from a distribution $F_n(\bth)$ that depends on the
true parameter $\bth\in\R^m$.  
Suppose a deterministic mapping $g_n$ produces the PTO estimator
\begin{equation*}
    \hbth \equiv \hbth_n = g_n\left(\mathcal D_n(\bth)\right).
\end{equation*}
Conditional on the observed sample, we generate $B$ bootstrap datasets
\[ \mathcal D_n^{\ast 1},\dots,\mathcal D_n^{\ast B} \;\stackrel{\text{i.i.d.}}{\sim}\; F_n(\hbth), \]
where recall that $F_n(\cdot)$ is the parametric model used to generate the dataset from the parameter.\footnote{We do not specify how the quantities determined by the nuisance parameters are generated, such as the residual in a linear regression model. In practice, certain bootstrap approaches can generate the dataset without estimating the nuisance. See Example~\ref{exp:pricing-bootstrap} later in the section for an example.}
Each resample yields a bootstrap version of the PTO estimator
$\hbth^{\ast b}=g_n(\mathcal D_n^{\ast b})$.
Given a candidate adjustment matrix $\Lambda$ we evaluate its performance in the
bootstrap world by
\begin{equation*}
    R^{\ast}_b(\Lambda)\coloneqq R_{\hbth}\left(\left(\I+\tfrac{\Lambda}{n}\right)\,\hbth^{\ast b}\right),\qquad b=1,\dots,B.
\end{equation*}
That is, the objective when the true parameter is $\hbth$ and the PTO estimator is $\hbth^{\ast b}$.
The empirical bootstrap objective is
\[ \hat R_B(\Lambda) = \frac1B\sum_{b=1}^B R^{\ast}_b(\Lambda). \]
Our data‑driven adjustment $\hat\Lambda$ is the minimizer of $\hat R_B$ over a pre‑specified search region
$\Omega\subset\mathbb R^{m\times m}$ (diagonal matrices, for example).
We provide the details in Algorithm~\ref{alg:bootstrapLambda}.
\begin{algorithm}[t]
\caption{Bootstrap‑based Estimation of the Adjustment Matrix $\hat\Lambda$}
\label{alg:bootstrapLambda}
\begin{algorithmic}[1]
    \Require Observed data $\mathcal D_n$, PTO estimator $\hat\theta = g_n(\mathcal D_n)$, \
            search domain $\Omega$, number of bootstrap draws $B$
    \For{$b=1$ \textbf{to} $B$}
        \State Generate bootstrap sample $\data_n^{\ast b} \sim F_n(\hbth)$
        \State Compute bootstrap PTO estimator $\hbth^{\ast b} \gets g_n(\data_n^{\ast b})$
    \EndFor
    \Statex
    \State Define $\hat R_B(\Lambda) \gets \frac{1}{B}\sum_{b=1}^B R_{\hbth}\bigl((\I+\Lambda/n)\,\hbth^{\ast b}\bigr)$
    \State Compute $\hat\Lambda \gets \arg\min_{\Lambda\in\Omega} \hat R_B(\Lambda)$ \hfill \emph{(e.g.\ grid search or coordinate descent)}
    \State \Return Adjusted estimator $ (\I+\hat\Lambda/n)\,\hat\theta$
\end{algorithmic}
\end{algorithm}

We use Example~\ref{exp:pricing} to illustrate the process. 
Although the dimension of the parameter is $m=1$ in Example~\ref{exp:pricing} and Proposition~\ref{prop:unknown-beta} has already provided a theoretically justified approach to find a data-driven adjustment, the steps using bootstrap are similar and the illustration can help to understand the general case.
\begin{example}[Data-driven Adjustment for Example~\ref{exp:pricing} Using Wild Bootstrap]\label{exp:pricing-bootstrap}
    Recall that the dataset consists of price-demand pairs: $\data = \{p_i, d_i\}_{i=1}^n$ with $d_i = a-\theta p_i+\epsilon_i$.
    We first obtain the OLS estimator $\hth$.
    Then we generate bootstrap datasets with $\data^{\ast b} = \{p_i,d^{\ast b}_i\}_{i=1}^n$ where
    \begin{align*}
        d^{\ast b}_{i} &= a - \hth p_i + v^{\ast b}_{i}\hat\epsilon_i,
    \end{align*}
    $v^{\ast b}_{i}\sim N(0,1)$ is an i.i.d. standard normal random variable, and $\hat\epsilon_i = d_i-a+\hth p_i$ is the residual of the $i$-th data point in the original dataset $\data$.
    We obtain an OLS estimator for each bootstrap dataset, denoted as $\hth^{\ast b}$.
    By Proposition~\ref{prop:known-beta}, we can obtain an oracle adjustment $\lambda(\hth)$ for the bootstrap datasets. 
    We can find
    \begin{align*}
        \hat\lambda =& \argmax_{\lambda } \frac{1}{B} \sum_{b=1}^B R_{\hth}\bigl((\I+\lambda/n)\,\hth^{\ast b}\bigr) 
    \end{align*}
    based on \eqref{eq:pricing-objective}.
    Note that in this case, the optimal $\hat\lambda$ can be analytically solved because the objective is a quadratic function of $1/(1+\lambda/n)$.
    Finally, we use $(1+\hat\lambda/n)\hth$ as the adjusted PTO estimator.
\end{example}
\paragraph{Why Does Bootstrap Work?} 
To provide theoretical justification for bootstrap, we compare the adjustment matrix $\hat\Lambda$ obtained from bootstrap and the oracle adjustment derived in Proposition~\ref{prop:multi-oracle}. 
When $A(\bth)$ is negative definite, as shown in Proposition~\ref{prop:multi-oracle}, we have $\Lambda = -\diag(A(\bth)^{-1} b(\bth))/2$.
In the bootstrap framework, conditional on $\hbth$ and taking $B\to\infty$, Algorithm~\ref{alg:bootstrapLambda} yields an adjustment matrix in the form:
\begin{equation*}
    \hat\Lambda = -\frac12\diag(\hat A(\hbth)^{-1} \hat b(\hbth)),
\end{equation*}
where
\begin{equation*}
\begin{aligned}
    \hat A(\hbth)_{ij} =& H_{\hbth}(\hbth)_{ij}\hth_i \hth_j,\\
    \hat b(\hbth) =& (H_{\hbth}(\hbth)\odot (2\Sigma_{\hat\theta}+M) )\one_m.
\end{aligned}
\end{equation*}
That is, we simply replace the true parameter $\bth$ by the PTO estimator $\hbth$.
This is because the bootstrap datasets are generated using $\hbth$ instead of $\bth$, which is used for the original dataset.
If we have 
\begin{enumerate}
    \item $\lim_{n\to\infty}\|\hbth_n- \bth\|_\infty=0$ almost surely,
    \item $H_{\bth}(\bth)$ and $\Sigma_{\bth}$ are continuous in $\bth$ with respect to $\|\cdot\|_\infty$,
\end{enumerate}
then we have 
\begin{equation*}
\begin{aligned}
   &\lim_{n\to\infty} -\frac12\diag(\hat A(\hbth)^{-1} \hat b(\hbth))= -\frac12\diag(A(\bth)^{-1} b(\bth))
\end{aligned}
\end{equation*}
almost surely by the continuous mapping theorem.
This result confirms that the bootstrap adjustment matrix will converge to the oracle adjustment matrix as $n$ grows.


In the next section, we test the bootstrap adjustment extensively alongside other policies. 
We can see that it outperforms PTO in every experiment and outperforms the oracle adjustment in certain cases.

\section{Numerical Experiments on Data-driven Pricing}\label{sec:case_study_linear_pricing}
We illustrate the practical value of the proposed adjustments on three pricing scenarios with large-scale simulation: linear demand with unknown price sensitivity (Example~\ref{exp:pricing}), linear demand with two unknown parameters (Example~\ref{exp:linear-demand-two-par}), and log-linear demand pricing (Example~\ref{exp:log-linear}).
In all the experiments, we compare {four} policies:
\begin{description}
\item[PTO:] the DM plugs the estimator $\hat{\bm\theta}_n$ into the optimal pricing
      formula.
\item[Oracle:] the estimator is adjusted by a factor $(\I+\Lambda/n)$ from
      Proposition~\ref{prop:known-beta}
      (single-parameter) or Proposition \ref{prop:multi-oracle} (multi-parameter).
\item[Plug-in:] the estimator is adjusted by a factor
      from Proposition~\ref{prop:unknown-beta} that is estimated from the data.
      It is available only in the single-parameter case.
\item[Bootstrap:] the estimator is adjusted in a data-driven fashion using Algorithm~\ref{alg:bootstrapLambda}. 
\end{description}
For every policy and experiment, we report the \emph{relative performance} of PTO:
${R_{\bth}(\hbth)}/{R_{\bth}(\bth)}$,
estimated by $10^{5}$ Monte Carlo instances.
Note that the quantity is always between 0 and 1; 
as $\hbth\to\bth$, the quantity converges to 1.
Moreover, we report the improvement of other policies relative to PTO:
\[
   \frac{R_{\bth}(\pi(\hbth))-R_{\bth}(\hbth)}{R_{\bth}(\hbth)},
\]
averaged over the same Monte Carlo instances.
A positive value indicates that the policy $\pi(\cdot)$ achieves better performance than PTO.

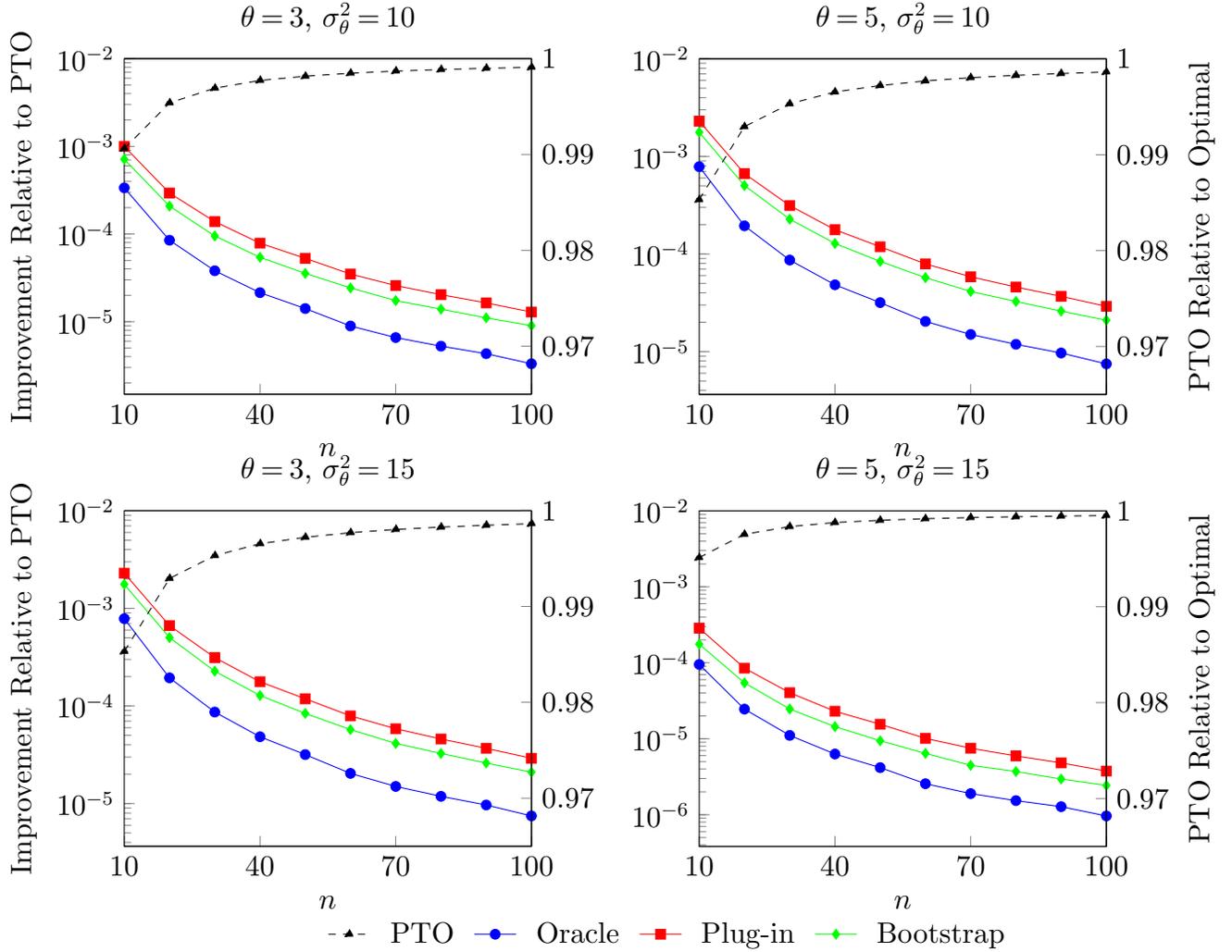
\begin{figure}[ht]
    \centering
    \begin{tikzpicture}
        \begin{scope}[xshift=0cm]                    
        \begin{semilogyaxis}[
                width=0.45\textwidth,
                xlabel={$n$}, xmin=10, xmax=100, 
                xtick={10,40,70,100}, xtick pos=left,
                ylabel={Improvement Relative to PTO},
                ymax = 0.01, ymin=0,
                ylabel near ticks,
                yticklabel pos=left,
                ytick pos=left,
                title={$\theta=3$, $\sigma_\theta^2=10$},
            ]
                \addplot [red, mark=square*] table[
                    x=n, y=dd, col sep=comma
                ] {./plot_data/theta3sigma10_2.csv};
                \label{plot:DataDriven}

                \addplot [blue, mark=*] table[
                    x=n,  y=oracle, col sep=comma,
                ] {./plot_data/theta3sigma10_2.csv};
                \label{plot:Oracle}

                \addplot [green, mark=diamond*] table[
                    x=n,  y=boot, col sep=comma,
                ] {./plot_data/theta3sigma10_2.csv};
                \label{plot:Boot}

            \end{semilogyaxis}
            \begin{axis}[
                width=0.45\textwidth,
                xmin=10, xmax=100, xmajorticks=false,
                ylabel={},
                ylabel near ticks,
                yticklabel pos=right,
                axis y line*=right,
                ymin=0.965, ymax=1, 
            ]
                \addplot [black, mark=triangle*, dashed] table[
                    x=n, y=pto, col sep=comma,
                ] {./plot_data/theta3sigma10_2.csv};
                \label{plot:PTO}
            \end{axis}
        \end{scope}

        \begin{scope}[xshift=0.5\textwidth]  
        \begin{semilogyaxis}[
                width=0.45\textwidth,
                xlabel={$n$}, xmin=10, xmax=100, 
                xtick={10,40,70,100}, xtick pos=left,
                ylabel={}, ymax=0.01, ymin=0,
                ylabel near ticks,
                yticklabel pos=left,
                ytick pos=left,
                title={$\theta=5$, $\sigma_\theta^2=10$},
            ]
                 \addplot [red, mark=square*] table[
                    x=n, y=dd, col sep=comma
                ] {./plot_data/theta3sigma15_2.csv};

                \addplot [blue, mark=*] table[
                    x=n,  y=oracle, col sep=comma,
                ] {./plot_data/theta3sigma15_2.csv};

                \addplot [green, mark=diamond*] table[
                    x=n,  y=boot, col sep=comma,
                ] {./plot_data/theta3sigma15_2.csv};
            \end{semilogyaxis}
            \begin{axis}[
                width=0.45\textwidth,
                xmin=10, xmax=100, xmajorticks=false,
                ylabel={PTO Relative to Optimal},
                ylabel near ticks,
                yticklabel pos=right,
                axis y line*=right,
                ymin=0.965, ymax=1, 
            ]
                \addplot [black, mark=triangle*, dashed] table[
                    x=n, y=pto, col sep=comma,
                ] {./plot_data/theta3sigma15_2.csv};
                \label{plot:PTO}
            \end{axis}
        \end{scope}

        \begin{scope}[yshift=-6.5cm]
            \begin{semilogyaxis}[
                width=0.45\textwidth,
                xlabel={$n$}, xmin=10, xmax=100, 
                xtick={10,40,70,100}, xtick pos=left,
                ylabel={Improvement Relative to PTO}, ymax=0.01, ymin=0,
                ylabel near ticks,
                yticklabel pos=left,
                ytick pos=left,
                title={$\theta=3$, $\sigma_\theta^2=15$},
            ]
                \addplot [red, mark=square*] table[
                    x=n, y=dd, col sep=comma
                ] {./plot_data/theta3sigma15_2.csv};

               \addplot [blue, mark=*] table[
                    x=n,  y=oracle, col sep=comma,
                ] {./plot_data/theta3sigma15_2.csv};

                \addplot [green, mark=diamond*] table[
                    x=n,  y=boot, col sep=comma,
                ] {./plot_data/theta3sigma15_2.csv};
                
            \end{semilogyaxis}
            \begin{axis}[
                 width=0.45\textwidth,
                xmin=10, xmax=100, xmajorticks=false,
                ylabel={},
                ylabel near ticks,
                yticklabel pos=right,
                axis y line*=right,
                ymin=0.965, ymax=1, 
            ]
               \addplot [black, mark=triangle*,dashed] table[
                    x=n,  y=pto, col sep=comma,
                ] {./plot_data/theta3sigma15_2.csv};
                \end{axis}
        \end{scope}

                \begin{scope}[xshift=0.5\textwidth, 
 yshift=-6.5cm]   
        \begin{semilogyaxis}[
                width=0.45\textwidth,
                xlabel={$n$}, xmin=10, xmax=100, 
                xtick={10,40,70,100}, xtick pos=left,
                ylabel={}, ymax=0.01, ymin=0,
                ylabel near ticks,
                yticklabel pos=left,
                ytick pos=left,
                title={$\theta=5$, $\sigma_\theta^2=15$},
            ]
                \addplot [red, mark=square*] table[
                    x=n, y=dd, col sep=comma
                ] {./plot_data/theta5sigma15_2.csv};

               \addplot [blue, mark=*] table[
                    x=n,  y=oracle, col sep=comma,
                ] {./plot_data/theta5sigma15_2.csv};

                \addplot [green, mark=diamond*] table[
                    x=n,  y=boot, col sep=comma,
                ] {./plot_data/theta5sigma15_2.csv};
            \end{semilogyaxis}
            \begin{axis}[
                width=0.45\textwidth,
                xmin=10, xmax=100, xmajorticks=false,
                ylabel={PTO Relative to Optimal},
                ylabel near ticks,
                yticklabel pos=right,
                axis y line*=right,
                ymin=0.965, ymax=1, 
            ]
                \addplot [black, mark=triangle*, dashed] table[
                    x=n, y=pto, col sep=comma,
                ] {./plot_data/theta5sigma15_2.csv};
                \label{plot:PTO}
            \end{axis}
        \end{scope}
    \end{tikzpicture}
   \begin{tikzpicture}[overlay, remember picture]
    \node[anchor=south,xshift=1.3cm, yshift=0cm] {
        \begin{tabular}{llllllll}
            \ref{plot:PTO} & PTO &
            \ref{plot:Oracle} & Oracle &
            \ref{plot:DataDriven} & Plug-in & 
            \ref{plot:Boot} & Bootstrap 
        \end{tabular}
    };
\end{tikzpicture}
    \caption{The performance against the sample size when the price sensitivity $\theta$ is unknown in the demand model $d(p)=a-\theta p$. 
    For PTO, we show the relative performance to the optimal revenue $\max_p\{p(a-\theta p)\}$ (right $y$-axis).
    For the other policies, their performance is computed relative to PTO and we plot the relative improvement in a log scale (left $y$-axis).}
    \label{fig:single-parameter}
\end{figure}

\textbf{Example~\ref{exp:pricing}: Linear Demand Pricing with Unknown Price Sensitivity}.
We consider the demand function $d(p)=a-\theta p$ where only $\theta$ is unknown.
We choose $a=60$ and $\theta\in\{3,5\}$ in the experiment.
For the historical data to estimate $\hth$, we generate $\data=\{p_i,d_i\}_{i=1}^n$ for $n\in\{10,20,\ldots,100\}$.
The price grid $p_i$ and realized demand $d_i$ is generated as follows:
\begin{enumerate}
    \item Given the data size $n$, we use a uniform grid of $n$ prices in the range $[0.1,6]$.
    \item For each price, we sample $d_i$ based on $d_i=a - \theta p_i + \epsilon_i$, where $\epsilon_i \sim N(0, \sigma_{\epsilon}^2)$.
    We test two values $\sigma_{\epsilon}^2 \in \{10, 15\}$ in the experiment.
    \item For each experimental setup $(\theta, \sigma_\epsilon, n)$, we simulate $10^5$ instances of data $\data$. This allows us to evaluate the average performance of our approaches and the PTO approach.
\end{enumerate}
For each simulated dataset $\data$, we estimate $\hth$ using OLS \eqref{eq:ols-pricing}.
The policies to be compared are
\begin{itemize}
    \item \emph{PTO}: we plug in the OLS estimator and use $a/(2\hth)$ as the price. 
    \item \emph{PTO with oracle adjustment}: following Proposition~\ref{prop:known-beta} in Section~\ref{sec:oracle-adj}, we first calculate the value of $\lambda=-\tfrac{(C+2)\sigma^2_\theta}{2\theta^2}$ in Proposition~\ref{prop:known-beta}. 
    Based on the form of $R_\theta(\hth)$ in \eqref{eq:pricing-objective}, we have $C=-6$. For $\sigma_\theta$ (the asymptotic standard deviation of $\hth$), the OLS formula leads to the expression $\sigma_\theta^2= \sigma_\epsilon^2n/(\sum_{i=1}^n p_i^2)$. As a result, we can calculate the value of $\lambda$ and the recommended price is $\hth(1+\lambda/n)$.
    \item \emph{PTO with data-driven plug-in adjustment}: following Proposition~\ref{prop:unknown-beta} in Section~\ref{sec:dd-adj}, we calculate the value of $\lambda=\tfrac{(2-C)\hsig^2_\theta}{2\hth_n^2}$ in Proposition~\ref{prop:unknown-beta}. The value $C=-6$ is the same as the oracle adjustment. For $\hth_n$ in the denominator, we simply use the OLS estimator. 
    For $\hsig_\theta^2$, we estimate it by $\hat{\sigma}_\epsilon^2n/(\sum_{i=1}^n p_i^2)$, where $\hat{\sigma}_\epsilon^2$ is computed using the mean squared error . 
    This allows us to calculate the value of $\lambda$ given the data and we then use the price $\hth(1+\lambda/n)$.
    \item \emph{PTO with bootstrap adjustment}:
    we follow Algorithm~\ref{alg:bootstrapLambda} with $B=10n$ and wild bootstrap as explained in Example~\ref{exp:pricing-bootstrap}.
    To construct the interval for the grid search, we use $\Omega=[-5\hat{\lambda},5\hat{\lambda}]$ with a grid size $0.1\hat\lambda$, where $\hat\lambda$ is the adjustment calculated in the plug-in method in the previous bullet point. 
    We select the adjustment that maximizes the average revenue in the bootstrap samples according to Step 5 and 6.
\end{itemize}
To measure the performance of the approaches above, note that the realized reward function $R_\theta(\pi(\hth))$ is $ p (a-\theta p)^+$,
where $p$ is the optimal price based on $\pi(\hth)$.
The optimal reward is given by $R_\theta(\theta) = a^2/(4\theta)$.
In Figure~\ref{fig:single-parameter}, we illustrate the performance of the four policies. 
The four panels show the setups with $\theta\in\{3,5\}$ and $\sigma_\epsilon^2\in\{10, 15\}$.

In each panel, we can observe a few consistent patterns.
First, across all data sizes, PTO with oracle, plug-in or bootstrap adjustments outperforms the standard PTO, i.e., their relative performance is above 0. 
This is consistent with our theoretical results in Propositions~\ref{prop:known-beta},~\ref{prop:unknown-beta} and Section~\ref{sec:bootstrap}. 
Second, PTO with data-driven adjustment (plug-in and bootstrap) outperforms that with the oracle adjustment.
Note that this doesn't hold in all applications.
Rather, it depends on the value of $C$. With $C=-6$, our theoretical results state that the data-driven adjustment will be better than the oracle adjustment (comparing \eqref{eq:single-par-improvement-gap} and \eqref{eq:single-par-dd-improvement}).
Third, although our theoretical framework is focused on the improvement when $n\to\infty$, it is more pronounced when the data size $n$ is small.

\textbf{Example~\ref{exp:linear-demand-two-par}: Linear Demand Pricing with Two Unknown Parameters}. 
Next, we explore the case when both the intercept and the price sensitivity are unknown in the demand function.
The dataset is generated in the same way as in the last example.
We use the OLS estimator $\hat{\theta}_1$ and $\hat{\theta}_2$.
For PTO, we simply plug in the estimated parameter and use the price ${\hat{\theta}_1}/(2\hat{\theta}_2)$.
For the oracle adjustment, we follow the discussion after Proposition~\ref{prop:multi-oracle} and in particular~\eqref{eq:foc-linear}.

For the bootstrap adjustment, we follow Algorithm~\ref{alg:bootstrapLambda} with\(B = 10n\) and bootstrap.
Moreover, we first obtain 
an initial estimate of the adjustment by a plug-in version of the oracle adjustment.
That is
\[
\hat{\lambda}_1 = - \frac{\hat{\Sigma}_{\theta,11}}{\hat{\theta}_1^2} + \frac{3\hat{\Sigma}_{\theta,12}}{\hat{\theta}_1\hat{\theta}_2} - \frac{2\hat{\Sigma}_{\theta,22}}{\hat{\theta}_2^2}, \quad \hat{\lambda}_2 = 0,
\]
where $\hat\Sigma$ is the estimated covariance matrix of $\hth_1$ and $\hth_2$, normalized by $n$, from the OLS estimator.
Starting from $\hat\lambda_1$ and $\hat\lambda_2$, we iteratively perform the coordinate descent with the range $[-5\hat\lambda_1,5\hat\lambda_1]\times [-0.5,0.5]$, with the grid size $0.1\hat\lambda_1$ and $0.01$.
We use the coordinate descent instead of the standard grid search to demonstrate the applicability when the parameter $\theta$ is high-dimensional.

\begin{figure}[htbp]
    \centering
    \begin{tikzpicture}
        \begin{scope}[xshift=0cm]                    
        \begin{semilogyaxis}[
                width=0.45\textwidth,
                xlabel={$n$}, xmin=10, xmax=100, 
                xtick={10,40,70,100}, xtick pos=left,
                ylabel={Improvement Relative to PTO},
                ymax = 0.01, ymin=0,
                ylabel near ticks,
                yticklabel pos=left,
                ytick pos=left,
                title={$\theta=3$, $\sigma_\theta^2=10$},
            ]

                \addplot [blue, mark=*] table[
                    x=n,  y=oracle, col sep=comma,
                ] {./plot_data/multi_theta3sigma10_2.csv};
                \label{plot:Oracle}
                 \addplot [green, mark=diamond*] table[
                    x=n,  y=boot, col sep=comma,
                ] {./plot_data/multi_theta3sigma10_2.csv};
                \label{plot:Boot}
            \end{semilogyaxis}
            \begin{axis}[
                width=0.45\textwidth,
                xmin=10, xmax=100, xmajorticks=false,
                ylabel={},
                ylabel near ticks,
                yticklabel pos=right,
                axis y line*=right,
                ymin=0.94, ymax=1, 
            ]
                \addplot [black, mark=triangle*, dashed] table[
                    x=n, y=pto, col sep=comma,
                ] {./plot_data/multi_theta3sigma10_2.csv};
                \label{plot:PTO}
            \end{axis}
        \end{scope}

        \begin{scope}[xshift=0.5\textwidth]  
        \begin{semilogyaxis}[
                width=0.45\textwidth,
                xlabel={$n$}, xmin=10, xmax=100, 
                xtick={10,40,70,100}, xtick pos=left,
                ylabel={}, ymax=0.01, ymin=0,
                ylabel near ticks,
                yticklabel pos=left,
                ytick pos=left,
                title={$\theta=5$, $\sigma_\theta^2=10$},
            ]
               \addplot [blue, mark=*] table[
                    x=n,  y=oracle, col sep=comma,
                ] {./plot_data/multi_theta5sigma10_2.csv};
               
                 \addplot [green, mark=diamond*] table[
                    x=n,  y=boot, col sep=comma,
                ] {./plot_data/multi_theta5sigma10_2.csv};
            \end{semilogyaxis}
            \begin{axis}[
                width=0.45\textwidth,
                xmin=10, xmax=100, xmajorticks=false,
                ylabel={PTO Relative to Optimal},
                ylabel near ticks,
                yticklabel pos=right,
                axis y line*=right,
                ymin=0.94, ymax=1, 
            ]
                \addplot [black, mark=triangle*, dashed] table[
                    x=n, y=pto, col sep=comma,
                ] {./plot_data/multi_theta5sigma10_2.csv};
                \label{plot:PTO}
            \end{axis}
        \end{scope}

        \begin{scope}[yshift=-6.5cm]
            \begin{semilogyaxis}[
               width=0.45\textwidth,
                xlabel={$n$}, xmin=10, xmax=100, 
                xtick={10,40,70,100}, xtick pos=left,
                ylabel={Improvement Relative to PTO},
                ylabel near ticks,
                yticklabel pos=left,
                title={$\theta=3$, $\sigma_\theta^2=15$},
            ]
               \addplot [blue, mark=*] table[
                    x=n,  y=oracle, col sep=comma,
                ] {./plot_data/multi_theta3sigma15_2.csv};
               
                 \addplot [green, mark=diamond*] table[
                    x=n,  y=boot, col sep=comma,
                ] {./plot_data/multi_theta3sigma15_2.csv};
            \end{semilogyaxis}
            \begin{axis}[
                width=0.45\textwidth,
                xmin=10, xmax=100, xmajorticks=false,
                ylabel={},
                ylabel near ticks,
                yticklabel pos=right,
                axis y line*=right,
                ymin=0.94, ymax=1,
            ]
               \addplot [black, mark=triangle*,dashed] table[
                    x=n,  y=pto, col sep=comma,
                ] {./plot_data/multi_theta3sigma15_2.csv};
                \end{axis}
        \end{scope}

                \begin{scope}[xshift=0.5\textwidth, 
 yshift=-6.5cm]   
        \begin{semilogyaxis}[
                width=0.45\textwidth,
                xlabel={$n$}, xmin=10, xmax=100, 
                xtick={10,40,70,100}, xtick pos=left,
                ylabel={}, ymax=0.01, ymin=0,
                ylabel near ticks,
                yticklabel pos=left,
                ytick pos=left,
                title={$\theta=5$, $\sigma_\theta^2=15$},
            ]
               \addplot [blue, mark=*] table[
                    x=n,  y=oracle, col sep=comma,
                ] {./plot_data/multi_theta5sigma15_2.csv};
               
                 \addplot [green, mark=diamond*] table[
                    x=n,  y=boot, col sep=comma,
                ] {./plot_data/multi_theta5sigma15_2.csv};
            \end{semilogyaxis}
            \begin{axis}[
                width=0.45\textwidth,
                xmin=10, xmax=100, xmajorticks=false,
                ylabel={PTO Relative to Optimal},
                ylabel near ticks,
                yticklabel pos=right,
                axis y line*=right,
                ymin=0.94, ymax=1, 
            ]
                \addplot [black, mark=triangle*, dashed] table[
                    x=n, y=pto, col sep=comma,
                ] {./plot_data/multi_theta5sigma15_2.csv};
                \label{plot:PTO}
            \end{axis}
        \end{scope}
    \end{tikzpicture}
   \begin{tikzpicture}[overlay, remember picture]
    \node[anchor=south, xshift=1.1cm, yshift=0cm] {
        \begin{tabular}{llllllll}
            \ref{plot:PTO} & PTO &
            \ref{plot:Oracle} & Oracle &
            \ref{plot:Boot} & Bootstrap 
        \end{tabular}
    };
\end{tikzpicture}
\caption{The comparison between PTO and the oracle/data-driven adjustment when both $a$ and $b$ are unknown in the demand function $d(p)=a-bp$. The PTO performance is relative to the optimal revenue $a^2/(4b)$ (the right $y$-axis). The performance of the adjustments is relative to PTO (the left $y$-axis). }
\label{fig:multi-linear-pricing}
\end{figure}
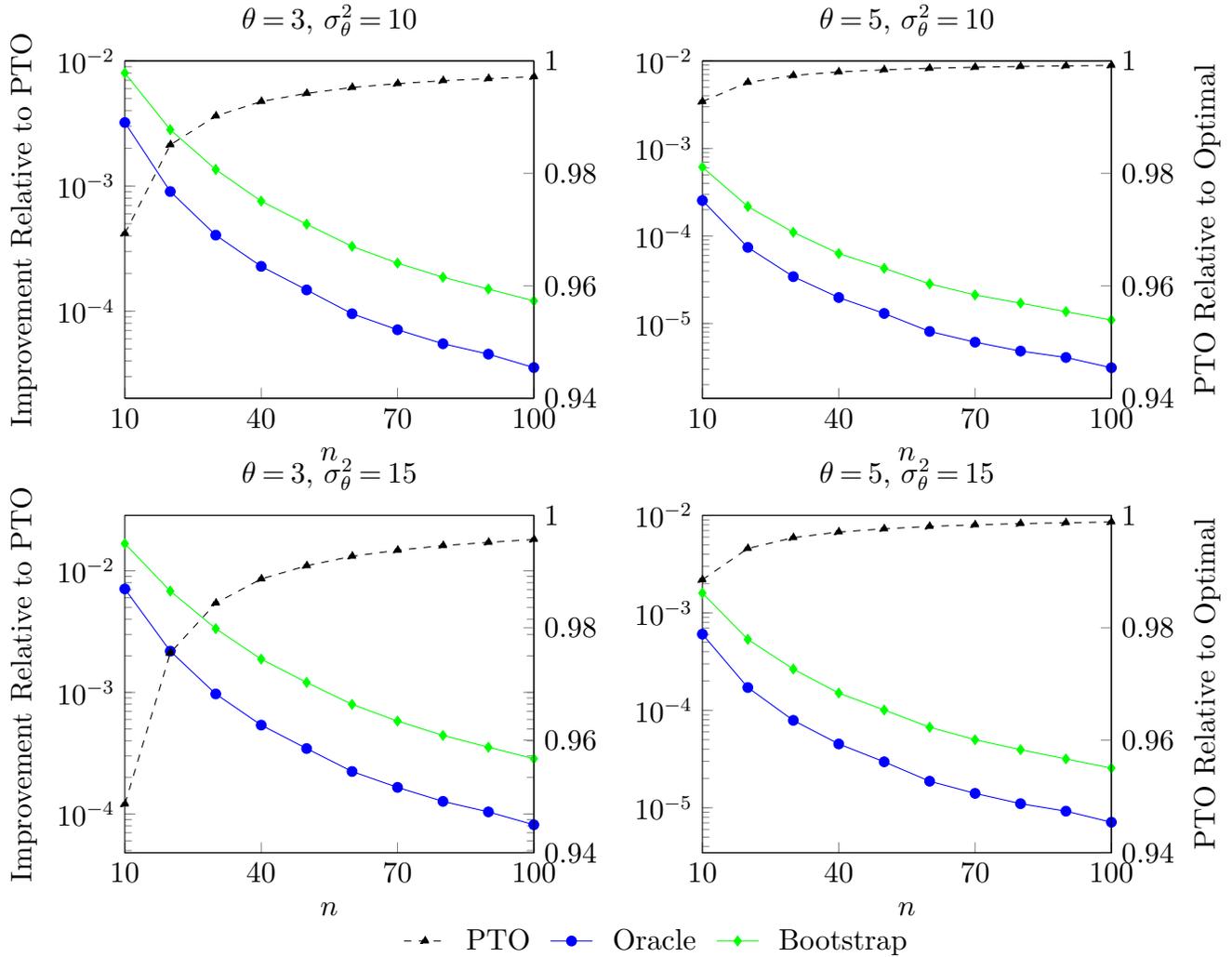

The comparison is shown in Figure~\ref{fig:multi-linear-pricing}.
The result is largely consistent: we observe that both adjustment schemes improve over PTO across all instances and data sizes, and the improvement is more significant for small datasets.
Overall, it demonstrates the robust performance of the adjustment schemes against the standard PTO approach.

\textbf{Example~\ref{exp:log-linear}: Log-linear Demand}.
In addition to the linear demand model in Example~\ref{exp:pricing}, we apply our framework to the log-linear demand model, $d(p) = \exp(a - \theta p)$, in Example~\ref{exp:log-linear}, where only the price sensitivity $\theta$ is unknown. 

We conduct a similar experiment to the linear demand.
We let $a = 8$ and $\theta \in \{3, 5\}$. 
Given sample size $n$, prices are chosen from a uniform grid over $[0.05,1]$.
The price range is smaller because the optimal price is $p^\ast = 1/\theta$ is small relative to the linear demand. 
We generate the historical demand $d_i$ 
using $\log(d_i) = a - \theta p_i + \epsilon_i$ for each price $p_i$, where $\epsilon_i$ is a normal random variable with a variance $\sigma^2_{\epsilon_i} \in \{0.5,1\}$. 
For each setup $(\theta, \sigma_\epsilon, n)$, we generate $100,000$ datasets to evaluate the average performance.

Given a dataset $\data = \{d_i, p_i\}_{i=1}^n$, we first estimate $\theta$ using OLS after transforming the demand:
\begin{equation}\label{eq:ols-loglinear}
\hth = -\frac{\sum_{i} p_i(\log(d_i)-a)}{\sum_i p_i^2}.
\end{equation}
For PTO, we simply set the price as $1/\hth$. For the oracle adjustment, we use the price $1/(\hth(1+\lambda/n))$, where
$\lambda = -{(C+2)\sigma_\theta^2}/(2\theta^2)$.
Note that $C = -4$ from \eqref{eq:log-linear-c} in Example~\ref{exp:log-linear}. 
For the data-driven plug-in adjustment, we use the formula 
$\lambda = {(2-C)\hat{\sigma}_\theta^2}/(2\hth^2)$
from Section~\ref{sec:dd-adj}, with $\hth$ estimated from \eqref{eq:ols-loglinear}. 
We estimate $\hat{\sigma}_\theta^2$ using the estimator for $\sigma_\epsilon$ and the OLS formula.
The recommended price is $1/(\hth(1 + \lambda/n))$.

\begin{figure}[htbp]
    \centering
    \begin{tikzpicture}
        \begin{scope}[xshift=0cm]                    
        \begin{semilogyaxis}[
                width=0.45\textwidth,
                xlabel={$n$}, xmin=10, xmax=100, 
                xtick={10,40,70,100}, xtick pos=left,
                ylabel={Improvement Relative to PTO},
                ylabel near ticks,
                yticklabel pos=left,
                ytick pos=left,
                title={$\theta=3$, $\sigma_\theta^2=0.5$},
            ]
                \addplot [red, mark=square*] table[
                    x=n, y=dd, col sep=comma
                ] {./plot_data/theta3sigma0.5_2.csv};
                \label{plot:DataDriven}

                \addplot [blue, mark=*] table[
                    x=n,  y=oracle, col sep=comma,
                ] {./plot_data/theta3sigma0.5_2.csv};
                \label{plot:Oracle}

                 \addplot [green, mark=diamond*] table[
                    x=n, y=boot, col sep=comma
                ] {./plot_data/theta3sigma0.5_2.csv};
                \label{plot:Boot}
                
            \end{semilogyaxis}
            \begin{axis}[
                width=0.45\textwidth,
                xmin=10, xmax=100, xmajorticks=false,
                ylabel={},
                ylabel near ticks,
                yticklabel pos=right,
                axis y line*=right,
               ymin=0.97, ymax=1, 
            ]
                \addplot [black, mark=triangle*, dashed] table[
                    x=n, y=pto, col sep=comma,
                ] {./plot_data/theta3sigma0.5_2.csv};
                \label{plot:PTO}
            \end{axis}
        \end{scope}

        \begin{scope}[xshift=0.5\textwidth]  
        \begin{semilogyaxis}[
                width=0.45\textwidth,
                xlabel={$n$}, xmin=10, xmax=100, 
                xtick={10,40,70,100}, xtick pos=left,
                ylabel={}, 
                ylabel near ticks,
                yticklabel pos=left,
                ytick pos=left,
                title={$\theta=5$, $\sigma_\theta^2=0.5$},
            ]
                 \addplot [red, mark=square*] table[
                    x=n, y=dd, col sep=comma
                ] {./plot_data/theta5sigma1_2.csv};

                \addplot [blue, mark=*] table[
                    x=n,  y=oracle, col sep=comma,
                ] {./plot_data/theta5sigma0.5_2.csv};

                 \addplot [green, mark=diamond*] table[
                    x=n, y=boot, col sep=comma
                ] {./plot_data/theta5sigma0.5_2.csv};
            \end{semilogyaxis}
            \begin{axis}[
                width=0.45\textwidth,
                xmin=10, xmax=100, xmajorticks=false,
                ylabel={PTO Relative to Optimal},
                ylabel near ticks,
                yticklabel pos=right,
                axis y line*=right,
                ymin=0.97, ymax=1, 
            ]
                \addplot [black, mark=triangle*, dashed] table[
                    x=n, y=pto, col sep=comma,
                ] {./plot_data/theta5sigma0.5_2.csv};
            \end{axis}
        \end{scope}

        \begin{scope}[yshift=-6.5cm]
            \begin{semilogyaxis}[
               width=0.45\textwidth,
                xlabel={$n$}, xmin=10, xmax=100, 
                xtick={10,40,70,100}, xtick pos=left,
                ylabel={Improvement Relative to PTO},
                ylabel near ticks,
                yticklabel pos=left,
                title={$\theta=3$, $\sigma_\theta^2=1$},
            ]
                \addplot [red, mark=square*] table[
                    x=n, y=dd, col sep=comma
                ] {./plot_data/theta3sigma1_2.csv};
               \addplot [blue, mark=*] table[
                    x=n,  y=oracle, col sep=comma,
                ] {./plot_data/theta3sigma1_2.csv};
                 \addplot [green, mark=diamond*] table[
                    x=n, y=boot, col sep=comma
                ] {./plot_data/theta3sigma1_2.csv};
                \end{semilogyaxis}
            \begin{axis}[
                width=0.45\textwidth,
                xmin=10, xmax=100, xmajorticks=false,
                ylabel={},
                ylabel near ticks,
                yticklabel pos=right,
                axis y line*=right,
                 ymin=0.97, ymax=1, 
            ]
               \addplot [black, mark=triangle*,dashed] table[
                    x=n,  y=pto, col sep=comma,
                ] {./plot_data/theta3sigma1_2.csv};
                \end{axis}
        \end{scope}

                \begin{scope}[xshift=0.5\textwidth, 
 yshift=-6.5cm]   
        \begin{semilogyaxis}[
                width=0.45\textwidth,
                xlabel={$n$}, xmin=10, xmax=100, 
                xtick={10,40,70,100}, xtick pos=left,
                ylabel={}, 
                ylabel near ticks,
                yticklabel pos=left,
                ytick pos=left,
                title={$\theta=5$, $\sigma_\theta^2=1$},
            ]
                \addplot [red, mark=square*] table[
                    x=n, y=dd, col sep=comma
                ] {./plot_data/theta5sigma1_2.csv};

                \addplot [blue, mark=*] table[
                    x=n,  y=oracle, col sep=comma,
                ] {./plot_data/theta5sigma1_2.csv};

                 \addplot [green, mark=diamond*] table[
                    x=n, y=boot, col sep=comma
                ] {./plot_data/theta5sigma1_2.csv};
            \end{semilogyaxis}
            \begin{axis}[
                width=0.45\textwidth,
                xmin=10, xmax=100, xmajorticks=false,
                ylabel={PTO Relative to Optimal},
                ylabel near ticks,
                yticklabel pos=right,
                axis y line*=right,
                 ymin=0.97, ymax=1, 
            ]
                \addplot [black, mark=triangle*, dashed] table[
                    x=n, y=pto, col sep=comma,
                ] {./plot_data/theta5sigma1_2.csv};
            \end{axis}
        \end{scope}
    \end{tikzpicture}
   \begin{tikzpicture}[overlay, remember picture]
    \node[anchor=south, xshift=1.3cm, yshift=0cm] {
       \begin{tabular}{llllllll}
            \ref{plot:PTO} & PTO &
            \ref{plot:Oracle} & Oracle &
            \ref{plot:DataDriven} & Plug-in & 
            \ref{plot:Boot} & Bootstrap 
        \end{tabular}
    };
\end{tikzpicture}
\caption{The comparison between PTO and the oracle/data-driven adjustment when $\theta$ is unknown in the demand function $d(p)=\exp(a-\theta p)$. The PTO performance is relative to the optimal revenue $\exp(a-1)/\theta$ (the right $y$-axis). The performance of the adjustment is relative to PTO (the left $y$-axis). }
\label{fig:loglinear}
\end{figure}
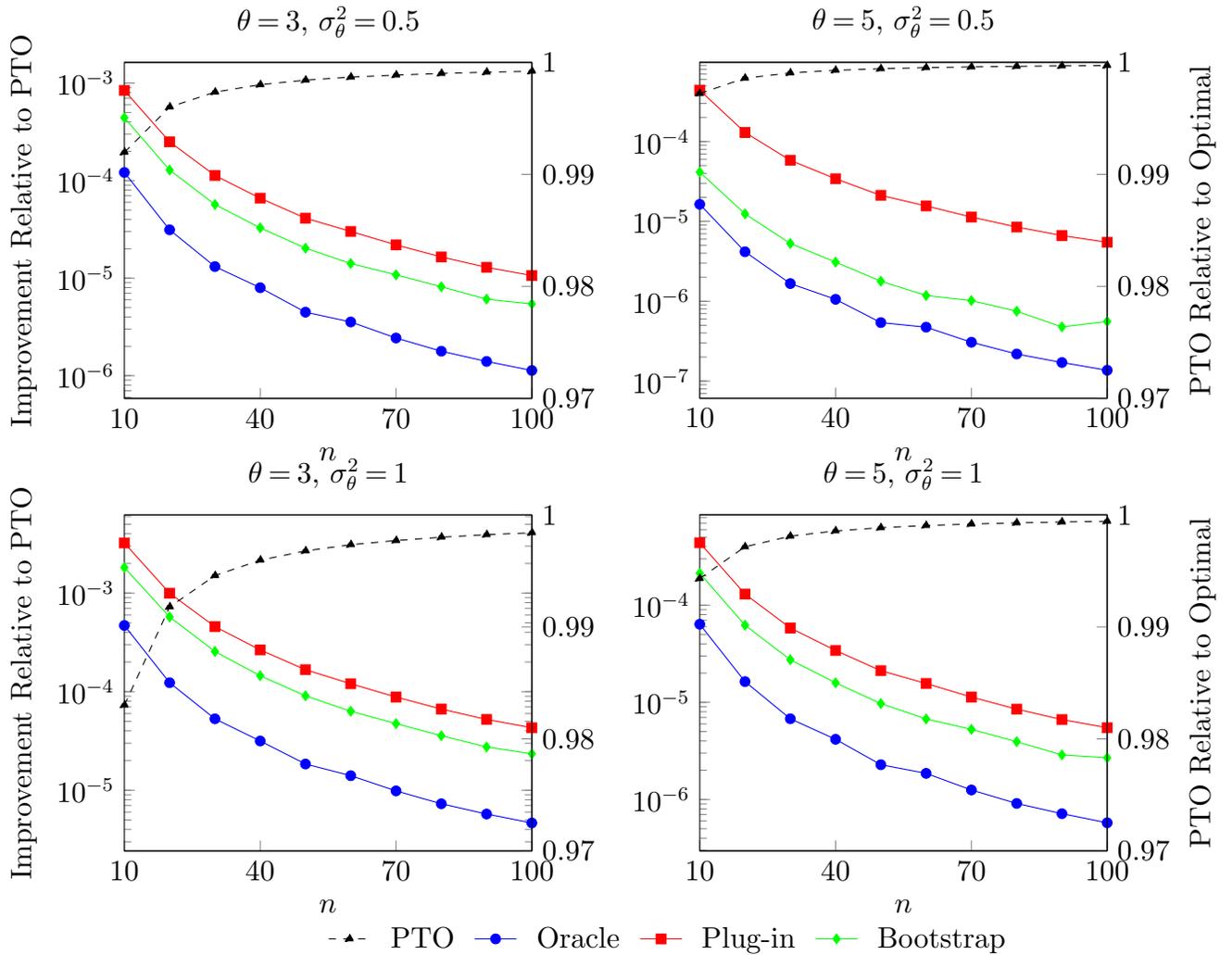

The bootstrap adjustment is selected by the similar procedure as the linear demand pricing.
The comparison of the performance of the four methods is shown in Figure~\ref{fig:loglinear}. 
We observe that both the oracle and data-driven adjustments outperform the PTO. This aligns with our findings in the linear demand case.

\section{Conclusion}

In this work, we propose a systematic modification to the standard predict-then-optimize (PTO) framework by incorporating a structured post-estimation adjustment that accounts for the curvature of the objective function. 
While our approach relies on specific conditions (most notably, Assumption~\ref{asp:r-derivative}) we establish that these conditions hold for a broad range of practical problems, particularly in pricing under various demand models.

We develop and analyze our adjustment across multiple settings: we derive oracle adjustments for single-parameter unbiased estimators, single-parameter biased estimators, and multi-parameter unbiased estimators. 
Additionally, we introduce a corresponding data-driven adjustment in the single-parameter unbiased setting. 
Interestingly, under certain conditions, our data-driven adjustment can be shown to outperform the oracle adjustment. 
Our numerical experiments confirm that this holds in both the linear and log-linear demand pricing problems, reinforcing our theoretical findings and demonstrating the effectiveness of our method across different data regimes.

An important avenue for future research is understanding the robustness of our adjustment to model misspecification. 
In real-world applications, the assumed demand model may not perfectly align with the true underlying consumer behavior. 
For instance, if the exponent $\gamma$ in the power-law demand model or the parameter $a$ in the single-parameter linear demand model is mis-specified, it remains unclear whether the proposed adjustment will still yield systematic improvement. 
A key question is whether our method provides guarantees that, for models sufficiently close to the true underlying decision problem, the adjustment continues to enhance decision quality. 
Investigating the sensitivity of our approach to small deviations in model parameters would offer deeper insights into its reliability and practical applicability, particularly in cases where firms must make pricing decisions with incomplete or imperfect information.


\bibliographystyle{informs2014} 
\bibliography{ref} 

%
%
%
\newpage
\begin{APPENDICES}
\section{Proofs in Section~\ref{sec:model}}\label{app:single-par-proof}
\begin{proof}{Proof of Proposition~\ref{prop:function-2-3derivative}:}
    In Proposition~\ref{prop:function-2-3derivative}, we first list of four objective functions, $R_\theta(\hth)$:
        \begin{itemize}
        \item $R_\theta(\hth)=f(\theta)\sum_{i=1}^I a_i \theta^{j_i}\hth^{K-j_i}+g(\theta) \hth+h(\theta)$ for $I,K\in \Z_+$, $j_i\in \R$ and an arbitrary function $f(\cdot)$. 
        \item $R_\theta(\hth)=f(\theta)\exp\left( \frac{a\hth}{\theta}\right)$ or $R_\theta(\hth)=\frac{f(\theta)}{\hth}\exp\left( \frac{a\theta}{\hth}\right)$ for $a\in \R$ and a function $f(\cdot)$.
         \item $R_\theta(\hth)=f(\theta)\log\left( g(\theta)\hth\right)+\hth$ for $g(\theta)\hth >0 $ and a function $f(\cdot)$. 
    \end{itemize}

The first function is $R_\theta(\hth)=f(\theta)\sum_{i=1}^I a_i \theta^{j_i}\hth^{K-j_i}+g(\theta) \hth+h(\theta)$ for $I,K\in \Z_+$, $j_i\in \R$ and an arbitrary function $f(\cdot)$. 
Taking the derivative with respect to $\hth$, we obtain the following derivative function:
\begin{equation*}
\begin{aligned}
R_\theta^{(1)}(\hth)=&f(\theta)\sum_{i=1}^I a_i \theta^{j_i}(K-j_i)\hth^{K-j_i-1}+g(\theta)\\
R_\theta^{(2)}(\hth)=&f(\theta)\sum_{i=1}^I a_i \theta^{j_i}(K-j_i)(K-j_i-1)\hth^{K-j_i-2}\\
R_\theta^{(3)}(\hth)=&f(\theta)\sum_{i=1}^I a_i \theta^{j_i}(K-j_i)(K-j_i-1)(K-j_i-2)\hth^{K-j_i-3}
\end{aligned}
\end{equation*}

Plugging $\hth=\theta$ into the derivative, Assumption~\ref{asp:r-derivative}(2) is satisfied by taking $C=\frac{\sum_{i=1}^I a_i (K-j_i)(K-j_i-1)(K-j_i-2)}{\sum_{i=1}^I a_i (K-j_i)(K-j_i-1)}$. Specifically, we have
\begin{equation*}
\begin{aligned}
\frac{R^{(3)}_{\theta}(\theta)}{R^{(2)}_{\theta}(\theta)}
=&\frac{\sum_{i=1}^I a_i (K-j_i)(K-j_i-1)(K-j_i-2)}{\sum_{i=1}^I a_i (K-j_i)(K-j_i-1)}\frac{1}{\theta}.
\end{aligned}
\end{equation*}

The second one is $R_\theta(\hth)=f(\theta)\exp\left( \frac{a\hth}{\theta}\right)$ for $a\in \R$ and a function $f(\cdot)$. Similarly, we obtained the second- and third-order derivative and evaluate at $\hth=\theta$,
\begin{equation*}
\begin{aligned}
R_\theta^{(2)}(\hth)=&\frac{a^2}{\theta^2}f(\theta)\exp\left( \frac{a\hth}{\theta}\right)\\
R_\theta^{(3)}(\hth)=&\frac{a^3}{\theta^3}f(\theta)\exp\left( \frac{a\hth}{\theta}\right)\\
\frac{R^{(3)}_{\theta}(\theta)}{R^{(2)}_{\theta}(\theta)}=&\frac{a}{\theta}.
\end{aligned}
\end{equation*}

The third one is $R_\theta(\hth)=\frac{f(\theta)}{\hth}\exp\left( \frac{a\theta}{\hth}\right)$ for $a\in \R$ and a function $f(\cdot)$. We have
\begin{equation*}
\begin{aligned}
\frac{R^{(3)}_{\theta}(\theta)}{R^{(2)}_{\theta}(\theta)}=&-\frac{(3a(a+2)+a(a^2+6a+6)+6a+6}{(a(a+2)+2a+2)}\frac{1}{\theta}
\end{aligned}
\end{equation*}

For $R_\theta(\hth)=f(\theta)\log\left( g(\theta)\hth\right)+\hth$ for $a\hth >0$ and a function $f(\cdot)$, 
\begin{equation*}
\begin{aligned}
R_\theta^{(2)}(\theta)=&-f(\theta)\frac{1}{\theta^2}\\
R_\theta^{(3)}(\theta)=&f(\theta)\frac{2}{\theta^3}\\
\frac{R^{(3)}_{\theta}(\theta)}{R^{(2)}_{\theta}(\theta)}=&\frac{-2}{\theta}
\end{aligned}
\end{equation*}

In the second part of Proposition~\ref{prop:function-2-3derivative}, the following objective functions $F(z, \theta)$ if $z(\hth)$ is derived from the first-order condition.
    \begin{itemize}
     \item $F(z, \theta) = a_1 \theta \log(z)-a_2z$, for $a_1\theta>0$ and $a_2>0$.
     \item $F(z, \theta) = f(a,\theta, \gamma)z(a-\theta z)^{\gamma}$, for known $\gamma$, $a>0$, $\theta>0$, and any $f(a,\theta, \gamma)$ can be any function independent of $z$. 
     \item $F(z, \theta) = f(a,\theta)\exp(a-\theta z)z$, for $a>0$, $\theta>0$, and any $f(a,\theta)$ can be any function independent of $z$.     
     \end{itemize}
     The first-order derivative of $F(z, \theta) = a_1 \theta \log(z)-a_2z$ is \(\frac{a_1\theta}{z}-a_2\). By plugging $z(\hth)=\frac{a_1}{a_2}\hth$ into our objective function, we obtain 
 $R_\theta(\hth)=a_1 \theta\log\left(\frac{a_1}{a_2}\hth\right)+a_1\hth$, which is a special case of $R_\theta(\hth)=f(\theta)\log\left( a\hth\right)+\hth$. 

The first-order derivative of $F(z, \theta) = f(a,\theta,\gamma)z(a-\theta z)^{\gamma}$ is 
$f(a,\theta,\gamma)\left((a-\theta z)^{\gamma}-\gamma \theta z(a-\theta z)^{\gamma-1} \right)$.
PTO decision is $z(\hth)=\frac{a}{\hth(1+\gamma)}$ and the objective function is
\[R_\theta(\hat{\theta})  
=f(a,\theta,\gamma)\left[\frac{a^{\gamma+1}}{(1+\gamma)\hat{\theta}}\left(1-\frac{\theta}{\hat{\theta}(1+\gamma)}\right)^\gamma\right].
\]
Taking the derivative with respect to $\hth$, we obtain the second-order derivative and the third-order derivative function,
\[
R^{(2)}_\theta(\theta)=f(a,\theta,\gamma)\left[-\frac{a^{\gamma+1}}{\gamma}\left(\frac{\gamma}{\gamma+1}\right)^{\gamma}\frac{1}{\theta^3}\right],\quad
R^{(3)}_\theta(\theta)=f(a,\theta,\gamma)\left[-\frac{a^{\gamma+1}}{\gamma^2}\left(\frac{\gamma}{\gamma+1}\right)^{\gamma}\frac{2(1+2\gamma)}{\theta^4}\right].
\]
Therefore, $\frac{R_\theta^{(3)}(\theta)}{R_\theta^{(2)}(\theta)}=\frac{C}{\theta}$ is given by
\[
\frac{R_\theta^{(3)}(\theta)}{R_\theta^{(2)}(\theta)}
= -\frac{2(1+2\gamma)}{\gamma}\frac{1}{\theta}.
\]

For $a>0$, $\theta>0$, and any $f(a,\theta)$ can be any function independent of $z$, the first-order derivative of $F(z, \theta) = f(a,\theta)\exp(a-\theta z)z$ is $f(a,\theta)\left(\exp(a-\theta z)-\theta z \exp(a-\theta z)\right)$. Hence, the plug-in decision is $z(\hth)=\frac{1}{\hth}$ and  the plug-in  objective function is 
\[
R_\theta(\hat{\theta}) = f(a,\theta)\left[exp(a-\frac{\theta}{\hth})\frac{1}{\hth}\right]
\]
$\frac{R_\theta^{(3)}(\theta)}{R_\theta^{(2)}(\theta)}=-\frac{4}{\theta}$, because the derivatives of $R_\theta(\hat{\theta})$ 
\[
R^{(2)}_\theta(\theta)=f(a,\theta)\left[-\frac{\exp(a-1)}{\theta^3}\right],\quad
R^{(3)}_\theta(\theta)=f(a,\theta)\left[\frac{4\exp(a-1)}{\theta^4}\right].
\]
\Halmos
\end{proof}

\begin{proof}{Proof for Proposition~\ref{prop:known-beta}:}

We first apply the fourth-order Taylor expansion of $R_{\theta}(\hth_n(1+\lambda^2/n))$ at $\hth_n(1+\lambda^2/n)=\theta$
with a fifth-order remainder:
\begin{equation*}
\begin{aligned}
R_{\theta}(\hth_n(1+\lambda/n))=&R(\theta)+\frac{R^{(2)}_{\theta}(\theta)}{2}(\hth_n+\hth_n\lambda/n-\theta)^2+\frac{R^{(3)}_{\theta}(\theta)}{6}(\hth_n+\hth_n\lambda/n-\theta)^3\\+&\frac{R^{(4)}_{\theta}(\theta)}{24}(\hth_n+\hth_n\lambda/n-\theta)^4+\frac{R^{(5)}_{\theta}(\theta')}{120}(\hth_n+\hth_n\lambda/n-\theta)^5.
\end{aligned}
\end{equation*}
Note that in the right-hand side, the random terms are $\hth_n$ and $\theta'$.
Taking expectation of the surrogate reward yields:
\begin{equation}\label{eq:generalR}
\begin{aligned}
\E[R_{\theta}(\hth_n(1+\lambda/n))]=&R(\theta)+\frac{R^{(2)}_{\theta}(\theta)}{2}\E[(\hth_n+\hth_n\lambda/n-\theta)^2]+\frac{R^{(3)}_{\theta}(\theta)}{6}\E[(\hth_n+\hth_n\lambda/n-\theta)^3]\\+&\frac{R^{(4)}_{\theta}(\theta)}{24}\E[(\hth_n+\hth_n\lambda/n-\theta)^4]+\E\left[\frac{R^{(5)}_{\theta}(\theta')}{120}(\hth_n+\hth_n\lambda/n-\theta)^5\right].
\end{aligned}
\end{equation}
The expectations in the first three terms in \eqref{eq:generalR} can be expressed as:
\begin{equation*}
\begin{aligned}
\E[(\hth_n+\hth_n\lambda/n-\theta)^2]=& \E[(1+\lambda/n)^2(\hth_n-\theta)^2]+\lambda^2\theta^2/n^2+2\E[(1+\lambda/n)(\hth_n-\theta)\lambda\theta/n ]\\
=&\sigma_n^2+\lambda^2\theta^2/n^2+2\lambda\sigma_n^2/n+o(n^{-2})\\
\E[(\hth_n+\hth_n\lambda/n-\theta)^3]=&\E[(\hth_n-\theta)^3]+3\lambda\theta\sigma_n^2/n+o(n^{-2})\\
\E[(\hth_n+\hth_n\lambda/n-\theta)^4]=&\E[(\hth_n-\theta)^4]+o(n^{-2}).
\end{aligned}
\end{equation*}
In the three equations, we have simply expanded them and repeatedly used the fact from Assumption~\ref{asp:eps-dist} that $\E[\hth_n-\theta]=o(1/n)$ and $\E[(\hth_n-\theta)^2]=o(1)$.
Under Assumption~\ref{asp:r-derivative}(3), we have 
\begin{equation*}
\begin{aligned}
\E\left[\frac{R^{(5)}_{\theta}(\theta')}{120}(\hth_n+\hth_n\lambda/n-\theta)^5\right] \leq \frac{M}{120}\E\left[(\hth_n+\hth_n\lambda/n-\theta)^5\right]=o(n^{-2}).
\end{aligned}
\end{equation*}
Substituting the above expectations into \eqref{eq:generalR}, the surrogate reward function becomes
\begin{equation}\label{eq:generalR-2}
\begin{aligned}
\E[R_{\theta}(\hth_n(1+\lambda/n))]=&R(\theta)+\frac{R^{(2)}_{\theta}(\theta)}{2}\sigma_n^2+\frac{R^{(3)}_{\theta}(\theta)}{6}\E[(\hth_n-\theta)^3]+\frac{R^{(4)}_{\theta}(\theta)}{24}\E[(\hth_n-\theta)^4]\\+&\frac{R^{(2)}_{\theta}(\theta)}{2}(\lambda^2\theta^2/n^2+2\lambda\sigma_n^2/n)+\frac{R^{(3)}_{\theta}(\theta)}{2}\lambda\theta\sigma_n^2/n+o(n^{-2}).
\end{aligned}
\end{equation}
Note that when we use $\lambda=0$, the expression above gives us the PTO performance $\E[R_{\theta}(\hth_n)]$.
Using $\frac{R^{(3)}_{\theta}(\theta)}{R^{(2)}_{\theta}(\theta)}=\frac{C}{\theta}$ from Assumption~\ref{asp:r-derivative}(2) and $\lim_{n\to \infty} n\sigma_n^2 = \sigma_\theta^2$ from Assumption~\ref{asp:eps-dist}, we have
\begin{equation}\label{eq:generalR-known-gap}
\begin{aligned}
n^2\E[R_{\theta}(\hth_n(1+\lambda/n))]-\E[R_{\theta}(\hth_n)]=&\frac{R^{(2)}_{\theta}(\theta)}{2}\left(\lambda^2\theta^2+2\lambda\sigma_\theta^2+C\lambda\sigma_\theta^2\right)+o(n^{-2}).
\end{aligned}
\end{equation}
The optimal $\lambda^\ast$ is $\frac{-(2+C)\sigma_\theta^2}{2\theta^2}$. 
Thus, the improvement gap in \eqref{eq:single-par-improvement-gap} is given by substituting $\lambda^\ast$ into \eqref{eq:generalR-known-gap}.
\Halmos
\end{proof}

\begin{proof}{Proof for Proposition \ref{prop:unknown-beta}}
The proof proceeds in four main steps. First, we express the adjusted estimator and perform a Taylor series expansion of the expected reward function around the true parameter $\theta$. Second, we analyze the moments of the terms in this expansion, carefully tracking their orders of magnitude. Third, we combine these results to derive an asymptotic expression for the improvement in expected reward, which we show is a quadratic function of the adjustment coefficient. Finally, we optimize this coefficient to maximize the improvement and verify that it matches the expression given in the proposition.

Let the data-driven plug-in adjustment coefficient be $\lambda_n = {k\hsig_\theta^2}/{\hth_n^2}$ for some constant $k$ to be determined. The adjusted estimator is $\tilde{\theta}_n = \hth_n(1+\lambda_n/n)$. We analyze the expected reward $\E[R_{\theta}(\tilde{\theta}_n)]$ by expanding $R_{\theta}(\cdot)$ around the true parameter $\theta$. The deviation of the adjusted estimator from the truth is
\begin{equation*}
\tilde{\theta}_n - \theta = \hth_n\left(1+\frac{k\hsig_\theta^2}{n\hth_n^2}\right) - \theta = (\hth_n - \theta) + \frac{k\hsig_\theta^2}{n\hth_n}.
\end{equation*}
We perform a fourth-order Taylor expansion of $R_{\theta}(\tilde{\theta}_n)$ around $\theta$. By Assumption~\ref{asp:r-derivative}(2), $R^{(1)}_{\theta}(\theta)=0$. Taking the expectation yields:
\begin{align}
\E[R_{\theta}(\tilde{\theta}_n)] = R_{\theta}(\theta) &+ \frac{R^{(2)}_{\theta}(\theta)}{2}\E\left[(\tilde{\theta}_n - \theta)^2\right] + \frac{R^{(3)}_{\theta}(\theta)}{6}\E\left[(\tilde{\theta}_n - \theta)^3\right] \notag \\
&+ \frac{R^{(4)}_{\theta}(\theta)}{24}\E\left[(\tilde{\theta}_n - \theta)^4\right] + \E\left[\text{Remainder}\right], \label{eq:proof-dd-taylor-exp}
\end{align}
where the remainder term is of order $o(n^{-2})$, as we will show. Our goal is to analyze each expectation on the right-hand side up to order $O(n^{-2})$.

Let $\delta_n = \hth_n - \theta$ and $\Delta_n = \frac{k\hsig_\theta^2}{n\hth_n}$. We need to compute $\E[(\delta_n + \Delta_n)^p]$ for $p=2,3,4$. We will repeatedly use Assumptions~\ref{asp:eps-dist}, \ref{asp:sigma-estimator}, and \ref{asp:betahat-away-from-zero}.

\emph{Analysis of the second moment ($p=2$):}
\begin{equation*}
\E[(\delta_n + \Delta_n)^2] = \E[\delta_n^2] + 2\E[\delta_n \Delta_n] + \E[\Delta_n^2].
\end{equation*}
The first term is the mean squared error of the original estimator, $\E[\delta_n^2] = \sigma_n^2$, which always appears regardless of the choice of $\lambda_n$.
For the cross-term, we analyze $\E[\delta_n \Delta_n] = \frac{k}{n}\E\left[(\hth_n-\theta)\frac{\hsig_\theta^2}{\hth_n}\right]$. We expand the term $\frac{\hsig_\theta^2}{\hth_n}$ around $(\theta, \sigma_\theta^2)$:
\begin{equation*}
\frac{\hsig_\theta^2}{\hth_n} = \frac{\sigma_\theta^2}{\theta} - \frac{\sigma_\theta^2}{\theta^2}(\hth_n-\theta) + \frac{1}{\theta}(\hsig_\theta^2-\sigma_\theta^2) + O((\hth_n-\theta)^2).
\end{equation*}
Multiplying by $(\hth_n-\theta)$ and taking expectations, we get:
\begin{align*}
\E\left[(\hth_n-\theta)\frac{\hsig_\theta^2}{\hth_n}\right] &= \frac{\sigma_\theta^2}{\theta}\E[\hth_n-\theta] - \frac{\sigma_\theta^2}{\theta^2}\E[(\hth_n-\theta)^2] + \frac{1}{\theta}\E[(\hth_n-\theta)(\hsig_\theta^2-\sigma_\theta^2)] + O(n^{-3/2}) \\
&= o(n^{-1}) - \frac{\sigma_\theta^2}{\theta^2}\left(\frac{\sigma_\theta^2}{n}\right) + \frac{1}{\theta}o(n^{-1}) + o(n^{-1}) = -\frac{\sigma_\theta^4}{n\theta^2} + o(n^{-1}),
\end{align*}
where we have used Assumptions~\ref{asp:eps-dist} and \ref{asp:sigma-estimator}. Thus, $2\E[\delta_n \Delta_n] = -\frac{2k\sigma_\theta^4}{n^2\theta^2} + o(n^{-2})$.
For the last term, $\E[\Delta_n^2] = \frac{k^2}{n^2}\E\left[\frac{\hsig_\theta^4}{\hth_n^2}\right]$. By the continuous mapping theorem and the law of large numbers, $\frac{\hsig_\theta^4}{\hth_n^2} \xrightarrow{p} \frac{\sigma_\theta^4}{\theta^2}$. With the moment conditions, we can show $\E\left[\frac{\hsig_\theta^4}{\hth_n^2}\right] = \frac{\sigma_\theta^4}{\theta^2} + o(1)$. Therefore, $\E[\Delta_n^2] = \frac{k^2\sigma_\theta^4}{n^2\theta^2} + o(n^{-2})$.
Combining these, the second moment is:
\begin{equation}
\E[(\tilde{\theta}_n - \theta)^2] = \E[(\hth_n-\theta)^2] - \frac{2k\sigma_\theta^4}{n^2\theta^2} + \frac{k^2\sigma_\theta^4}{n^2\theta^2} + o(n^{-2}). \label{eq:proof-dd-moment2}
\end{equation}

\emph{Analysis of the third moment ($p=3$):}
\begin{equation*}
\E[(\delta_n + \Delta_n)^3] = \E[\delta_n^3] + 3\E[\delta_n^2 \Delta_n] + 3\E[\delta_n \Delta_n^2] + \E[\Delta_n^3].
\end{equation*}
The leading term that depends on $k$ is $3\E[\delta_n^2 \Delta_n] = \frac{3k}{n}\E\left[(\hth_n-\theta)^2\frac{\hsig_\theta^2}{\hth_n}\right]$. A similar expansion shows:
\begin{equation*}
\E\left[(\hth_n-\theta)^2\frac{\hsig_\theta^2}{\hth_n}\right] = \E\left[(\hth_n-\theta)^2\left(\frac{\sigma_\theta^2}{\theta} + o_p(1)\right)\right] = \frac{\sigma_\theta^2}{\theta}\E[(\hth_n-\theta)^2] + o(n^{-1}) = \frac{\sigma_\theta^4}{n\theta} + o(n^{-1}).
\end{equation*}
The other terms are of smaller order: $\E[\delta_n \Delta_n^2] = O(n^{-5/2})$ and $\E[\Delta_n^3] = O(n^{-3})$. Thus:
\begin{equation}
\E[(\tilde{\theta}_n - \theta)^3] = \E[(\hth_n-\theta)^3] + \frac{3k\sigma_\theta^4}{n^2\theta} + o(n^{-2}). \label{eq:proof-dd-moment3}
\end{equation}

\emph{Analysis of higher moments ($p \ge 4$):}
All terms involving $\Delta_n$ in the expansion of $\E[(\delta_n + \Delta_n)^p]$ for $p \ge 4$ are of order $o(n^{-2})$. For instance, $\E[(\delta_n + \Delta_n)^4] = \E[\delta_n^4] + 4\E[\delta_n^3 \Delta_n] + \dots = \E[\delta_n^4] + O(n^{-5/2})$. Therefore,
\begin{align}
\E[(\tilde{\theta}_n - \theta)^4] &= \E[(\hth_n-\theta)^4] + o(n^{-2}), \label{eq:proof-dd-moment4} \\
\E\left[\text{Remainder}\right] &= O(\E[(\tilde{\theta}_n - \theta)^5]) = O(n^{-5/2}) = o(n^{-2}). \label{eq:proof-dd-moment5}
\end{align}

We now substitute the moment calculations \eqref{eq:proof-dd-moment2}-\eqref{eq:proof-dd-moment5} back into the Taylor expansion \eqref{eq:proof-dd-taylor-exp}. The terms not involving $k$ reconstruct the expected reward of the PTO estimator, $\E[R_\theta(\hth_n)]$, up to order $o(n^{-2})$. The difference in expected reward is therefore driven by the terms involving $k$:
\begin{align*}
\E[R_{\theta}(\tilde{\theta}_n)] - \E[R_{\theta}(\hth_n)] &= \frac{R^{(2)}_{\theta}(\theta)}{2}\left(-\frac{2k\sigma_\theta^4}{n^2\theta^2} + \frac{k^2\sigma_\theta^4}{n^2\theta^2}\right) + \frac{R^{(3)}_{\theta}(\theta)}{6}\left(\frac{3k\sigma_\theta^4}{n^2\theta}\right) + o(n^{-2}) \\
&= \frac{1}{n^2}\left[ \frac{R^{(2)}_{\theta}(\theta)\sigma_\theta^4}{2\theta^2}(k^2 - 2k) + \frac{R^{(3)}_{\theta}(\theta)\sigma_\theta^4}{2\theta}k \right] + o(n^{-2}).
\end{align*}
Now, we use the condition from Assumption~\ref{asp:r-derivative}(2), which states $\frac{R^{(3)}_{\theta}(\theta)}{R^{(2)}_{\theta}(\theta)}=\frac{C}{\theta}$. This allows us to write $R^{(3)}_{\theta}(\theta) = \frac{C}{\theta}R^{(2)}_{\theta}(\theta)$. Substituting this into the expression gives:
\begin{align*}
\E[R_{\theta}(\tilde{\theta}_n)] - \E[R_{\theta}(\hth_n)] &= \frac{1}{n^2}\left[ \frac{R^{(2)}_{\theta}(\theta)\sigma_\theta^4}{2\theta^2}(k^2 - 2k) + \frac{(C/\theta)R^{(2)}_{\theta}(\theta)\sigma_\theta^4}{2\theta}k \right] + o(n^{-2}) \\
&= \frac{R^{(2)}_{\theta}(\theta)\sigma_\theta^4}{2n^2\theta^2} \left[ (k^2 - 2k) + Ck \right] + o(n^{-2}) \\
&= \frac{R^{(2)}_{\theta}(\theta)\sigma_\theta^4}{2n^2\theta^2} \left[ k^2 + (C-2)k \right] + o(n^{-2}).
\end{align*}

To maximize the improvement, we must choose $k$ to optimize the quadratic term $ k^2 + (C-2)k$. Since Assumption~\ref{asp:r-derivative}(2) states $R^{(2)}_{\theta}(\theta) < 0$, we need to {minimize} it. 
The minimum is achieved at $k$ satisfying $2k + (C-2) = 0$, which gives the optimal coefficient:
\begin{equation*}
k^\ast = \frac{2-C}{2}.
\end{equation*}
This corresponds to the proposed adjustment coefficient $\lambda_n^\ast = \frac{(2-C)\hsig_\theta^2}{2\hth_n^2}$.
Substituting $k^\ast$ back into the improvement formula, the value of the quadratic at the optimum is:
\begin{equation*}
(k^\ast)^2 + (C-2)k^\ast = -\frac{(2-C)^2}{4}.
\end{equation*}
Therefore, the asymptotic improvement is:
\begin{equation*}
\lim_{n\to\infty} n^2 \left(\E\left[R_{\theta}\left(\hth_n\left(1+\frac{\lambda_n^\ast}{n}\right)\right)\right]-\E[R_{\theta}(\hth_n)]\right) = \frac{R^{(2)}_{\theta}(\theta)\sigma_\theta^4}{2\theta^2} \left( -\frac{(2-C)^2}{4} \right) = -\frac{R^{(2)}_{\theta}(\theta)(2-C)^2\sigma_\theta^4}{8\theta^2}.
\end{equation*}
Since $R^{(2)}_{\theta}(\theta) < 0$, this improvement is non-negative, which completes the proof.
\hfill\Halmos

\end{proof}

\section{Proofs in Section~\ref{sec:multi-para-oracle}}
\label{app:multi-par}
\begin{proof}{Proof of Proposition~\ref{prop:multi-oracle}:}
We will present the following easy-to-verify properties throughout the proof.
For vectors $a,b\in \R^m$, we have
\begin{equation}\label{eq:tensor-norm}
    \|a\otimes b\|_{\max} \le \|a\|_\infty \|b\|_\infty\\
\end{equation}
where $\|A\|_{\max}\coloneqq \max_{i_1,i_2,\dots}|A_{i_1i_2\dots}|$ for a tensor/matrix and  $\|a\|_\infty=\max_i |a_i|$ for a vector.
It is also easy to see that
\begin{equation}\label{eq:Ax-norm}
    \|Ax\|_{\infty} \le m \|A\|_{\max}\|x\|_\infty,
\end{equation}
where $A\in \R^{m\times m}$ and $x\in \R^m$.

We next apply the multivariate Taylor expansion of $\E[R_{\bm \theta}((\I+\Lambda/n)\hbth_n)]$ 
at $(\I+\Lambda/n)\hbth_n=\bm \theta$, up to a fifth-order remainder.
We will derive the expression for a general $\Lambda$ without using the property that it is a diagonal matrix.
\begin{equation}\label{eq:multi-taylor}
\begin{aligned}
    \E[R_{\bm \theta}((\I+ \Lambda/n)\hbth_n )]=&R(\bm \theta)+\frac{1}{2}\left\langle \Htheta,\E\left[\left((\I+\Lambda/n)\hbth_n-\bm \theta\right)^{\otimes 2}\right]\right\rangle \\
    +&\frac{1}{6}\left\langle \nabla^3 R_{\bm{\theta}}(\bm{\theta}), \E\left[\left((\I+\Lambda/n)\hbth_n-\bm \theta\right)^{\otimes 3}\right]\right\rangle\\
    +&\frac{1}{24}\left\langle \nabla^4 R_{\bm{\theta}}(\bm{\theta}), \E\left[\left((\I+\Lambda/n)\hbth_n-\bm \theta\right)^{\otimes 4}\right]\right\rangle \\
    +& \frac{1}{120} \E\left[\left\langle \nabla^5 R_{\bm{\theta}}(\bm{\theta}'), \left((\I+\Lambda/n)\hbth_n-\bm \theta\right)^{\otimes 5}\right\rangle\right],
    \end{aligned}
\end{equation}
where $\bm \theta'\in \R^m$ is a random vector that is between $\bm \theta$ and $\bm \hth_n$.
We further expand the terms in the right-hand side of \eqref{eq:multi-taylor} and focus on the terms that are not $o(n^{-2})$.
Note that we don't have the first-order term because of part (1) of Assumption~\ref{asp:r-multi-derivative}.

We first tackle the second-order term: 
\begin{align}
    \E\left[\left((\I+\Lambda/n)\hbth_n-\bm \theta\right)^{\otimes 2}\right] =& \E\left[\left((\hbth_n-\bm \theta)+\Lambda(\hbth_n-\bm \theta)/n+\Lambda \bm \theta/n\right)^{\otimes 2}\right] \notag\\
    =&\E\left[(\hbth_n-\bm \theta)^{\otimes 2}\right]+ \frac{1}{n^2}\E\left[(\Lambda(\hbth_n-\bm \theta))^{\otimes 2}\right]+\frac{1}{n^2}(\Lambda\bm \theta)^{\otimes 2}\notag\\
    & + \frac{1}{n}\E\left[(\hbth_n-\bm \theta)\otimes (\Lambda(\hbth_n-\bm \theta))\right]+\frac{1}{n}\E\left[\hbth_n-\bm \theta\right]\otimes (\Lambda\bm \theta)\notag\\
    & +\frac{1}{n}\E\left[\Lambda(\hbth_n-\bm \theta)\otimes (\hbth_n-\bm \theta)\right]+\frac{1}{n^2}\Lambda\E\left[\hbth_n-\bm \theta\right]\otimes (\Lambda\bm \theta)\notag\\
    & +\frac{1}{n}\Lambda\bm\theta\otimes \E\left[\hbth_n-\bm \theta\right]+\frac{1}{n^2}\Lambda\bm \theta\otimes \left(\Lambda\E\left[\hbth_n-\bm \theta\right]\right).\label{eq:multi-taylor-second}
\end{align}
In the right-hand side, the first term is $\Sigma_n$ according to the definition in part (2) of Assumption~\ref{asp:multi-eps-dist}.
For the second term, using \eqref{eq:Ax-norm} and then \eqref{eq:tensor-norm}, we have that 
\begin{align*}
    \left\| \E\left[(\Lambda(\hbth_n-\bm \theta))^{\otimes 2}\right]\right\|_{\max}\le \E\left[\|\Lambda(\hbth_n-\bm \theta)\|_\infty^2\right] \le m^2 \|\Lambda\|_{\max}^2\E[\|\hbth-\bm\theta\|_\infty^2].
\end{align*}
Again, $\E[\|\hbth-\bm\theta\|_\infty^2]=O(n^{-1})$ by part (2) of Assumption~\ref{asp:multi-eps-dist}.
Hence the second term of \eqref{eq:multi-taylor-second} is $o(n^{-2})$.
For the fourth term, we have 
\begin{align*}
    \E\left[(\hbth_n-\bm \theta)\otimes (\Lambda(\hbth_n-\bm \theta))\right] = \E\left[ (\hbth_n-\bm \theta)(\hbth_n-\bm \theta)^\top\Lambda^\top\right]=\Sigma_n\Lambda^\top=\frac{1}{n}\Sigma_\theta\Lambda^\top+o(n^{-1}),
\end{align*}
where we use part (2) of Assumption~\ref{asp:multi-eps-dist}.
The fifth term is $o(n^{-2})$ because of part (1) of Assumption~\ref{asp:multi-eps-dist} and property~\eqref{eq:tensor-norm}.
Similar to the fourth term, the sixth term is $\Lambda \Sigma_\theta/n^2+o(n^{-2})$.
It is easy to see that the other terms are $o(n^{-2})$.
Thus, we can simplify \eqref{eq:multi-taylor-second} to
\begin{align*}
    \eqref{eq:multi-taylor-second} = \Sigma_n+ \frac{1}{n^2}(\Lambda\bm\theta)^{\otimes 2} +\frac{1}{n^2}(\Lambda\Sigma_\theta+\Sigma_\theta\Lambda^\top) + o(n^{-2}).
\end{align*}
Plugging it into \eqref{eq:multi-taylor}, we have 
\begin{align}\label{eq:taylor-second-final}
    \left\langle \Htheta,\E\left[\left((\I+\Lambda/n)\hbth_n-\bm \theta\right)^{\otimes 2}\right]\right\rangle=& \left\langle H_\theta(\bm\theta), \Sigma_n\right\rangle+  \frac{2}{n^2}\left\langle H_\theta(\bm\theta), \Lambda\Sigma_\theta\right\rangle
    + \frac{1}{n^2}\left\langle H_\theta(\bm\theta), \Lambda \theta \theta^\top \Lambda^\top\right\rangle,
\end{align}
where we have used the property that $\langle H_\theta(\bm\theta), A\rangle = \langle H_\theta(\bm\theta), A^\top\rangle$ since the Hessian matrix is symmetric.

Next we use the same argument for the third-order expansion in \eqref{eq:multi-taylor}, and obtain
\begin{align}
    \E\left[\left((\I+\Lambda/n)\hbth_n-\bm \theta\right)^{\otimes 3}\right] = &\E[ (\hbth_n-\bm\theta)^{\otimes 3}] + \frac{1}{n}\E[(\hbth_n-\bm \theta)\otimes (\hbth_n-\bm \theta)\otimes (\Lambda\bm\theta)]\notag\\
    &+\frac{1}{n}\E[(\hbth_n-\bm \theta)\otimes  (\Lambda\bm\theta) \otimes(\hbth_n-\bm \theta)]+\frac{1}{n}\E[(\Lambda\bm\theta) \otimes(\hbth_n-\bm \theta)\otimes  (\hbth_n-\bm \theta)] \notag\\
    &+ o(n^{-2})    \notag\\
     =& \E[ (\hbth_n-\bm\theta)^{\otimes 3}] + \frac{1}{n}(\Sigma_n\otimes (\Lambda \bm\theta)+\Lambda\bm\theta\otimes\Sigma_n)\notag\\
    &+\frac{1}{n}\E[(\hbth_n-\bm \theta)\otimes  (\Lambda\bm\theta) \otimes(\hbth_n-\bm \theta)]+o(n^{-2})\label{eq:multi-taylor-third}
\end{align}
To plug it into \eqref{eq:multi-taylor}, we prove the following claim:
\begin{align}\label{eq:order3-permutation}
    \langle \nabla^3 R_{\bm{\theta}}(\bm{\theta}), \Sigma_n\otimes (\Lambda \bm\theta)\rangle = \langle \nabla^3 R_{\bm{\theta}}(\bm{\theta}), (\Lambda \bm\theta)\otimes \Sigma_n\rangle =  \langle \nabla^3 R_{\bm{\theta}}(\bm{\theta}), (\hbth_n-\bm \theta)\otimes  (\Lambda\bm\theta) \otimes(\hbth_n-\bm \theta)\rangle.
\end{align}
More generally, we prove that for $a_1,a_2,a_3\in \R^{m}$ and a permutation of the three vectors $b_1,b_2,b_3$, we have
\begin{equation}\label{eq:R3-permutation}
    \langle \nabla^3 R_{\bm{\theta}}(\bm{\theta}), a_1\otimes a_2\otimes a_3\rangle = 
    \langle \nabla^3 R_{\bm{\theta}}(\bm{\theta}), b_1\otimes b_2\otimes b_3\rangle.
\end{equation}
This is true because of the following property of $\nabla^3 R_{\bm{\theta}}(\bm{\theta})$:
for $i_1,i_2,i_3\in \{1,\dots,m\}$ and a permutation $j_1,j_2,j_3$, we have 
$\nabla^3 R_{\bm{\theta}}(\bm{\theta})_{i_1i_2i_3}=\nabla^3 R_{\bm{\theta}}(\bm{\theta})_{j_1j_2j_3}$.
We can direct verify \eqref{eq:R3-permutation} using this property. 
As a result, if $a_1=\hbth_n-\bm\theta$, $a_2=\Lambda\bm\theta$, and $a_3=\hbth_n-\bm\theta$, then \eqref{eq:R3-permutation} implies \eqref{eq:order3-permutation}.
Plugging \eqref{eq:multi-taylor-third} into \eqref{eq:multi-taylor}, we have
\begin{align}
    \left\langle \nabla^3 R_{\bm{\theta}}(\bm{\theta}), \E\left[\left((\I+\Lambda/n)\hbth_n-\bm \theta\right)^{\otimes 3}\right]\right\rangle 
    =&\left\langle \nabla^3 R_{\bm\theta}(\bm\theta),\E[(\hbth_n-\bm \theta)^{\otimes 3}]\right\rangle + \frac{3}{n}  \left\langle \nabla^3 R_{\bm{\theta}}(\bm{\theta}),\Lambda\bm\theta\otimes\Sigma_n\right\rangle+o(n^{-2})\notag\\
    =&\left\langle \nabla^3 R_{\bm\theta}(\bm\theta),\E[(\hbth_n-\bm \theta)^{\otimes 3}]\right\rangle + \frac{3}{n^2}  \left\langle \nabla^3 R_{\bm{\theta}}(\bm{\theta}),\Lambda\bm\theta\otimes\Sigma_\theta\right\rangle+o(n^{-2})\notag\\
    =&\left\langle \nabla^3 R_{\bm\theta}(\bm\theta),\E[(\hbth_n-\bm \theta)^{\otimes 3}]\right\rangle + \frac{3}{n^2}  \left\langle H_\theta(\bm\theta),\Lambda M\right\rangle+o(n^{-2}).\label{eq:taylor-third-final}
\end{align}
In the second inequality, we use part (2) of Assumption~\ref{asp:multi-eps-dist} which implies $\Sigma_n = \Sigma_\theta/n+o(n^{-1})$.
In the last step, we apply part (2) of Assumption~\ref{asp:r-multi-derivative}.

Next we analyze the fourth-order term in \eqref{eq:multi-taylor}.
There is only one term that is not $o(n^{-2})$, and we have
\begin{equation}\label{eq:taylor-four-final}
\left\langle \nabla^4 R_{\bm{\theta}}(\bm{\theta}), \E\left[\left((\I+\Lambda/n)\hbth_n-\bm \theta\right)^{\otimes 4}\right]\right\rangle 
    =\left\langle \nabla^4 R_{\bm\theta}(\bm\theta),\E[(\hbth_n-\bm \theta)^{\otimes 4}]\right\rangle+o(n^{-2})
\end{equation}
For the fifth-order term, because of part (3) of Assumption~\ref{asp:r-multi-derivative}, we have $\|\nabla^5 R_{\bm{\theta}}(\bm{\theta}')\|_{\max}\le K$ for some constant $K$.
As a result, we have
\begin{align}  
\E\left[\left\langle \nabla^5 R_{\bm{\theta}}(\bm{\theta}'), \left((\I+\Lambda/n)\hbth_n-\bm \theta\right)^{\otimes 5}\right\rangle\right]&\le K m^5 \left\|\E\left[\left((\I+\Lambda/n)(\bm \hth_n-\bm \theta)+\Lambda \bm \theta/n\right)^{\otimes 5}\right]\right\|_{\max}\notag\\
&= Km^5 \|\E[(\hbth_n-\bm\theta)^{\otimes 5}]\|_{\max}+o(n^{-2})=o(n^{-2}),\label{eq:taylor-five-final}
\end{align}
where the last equality is due to part (3) of Assumption~\ref{asp:multi-eps-dist}.

Combining the terms \eqref{eq:taylor-second-final}, \eqref{eq:taylor-third-final}, \eqref{eq:taylor-four-final}, and \eqref{eq:taylor-five-final} and rearranging the terms based on $\Lambda$, we have 
\begin{align*}
    \eqref{eq:multi-taylor} = &R(\bm \theta)+\frac{1}{2}\left\langle \Htheta, \Sigma_n\right\rangle+\frac{1}{6}\left\langle \nabla^3 R_{\bm\theta}(\bm\theta),\E[(\hbth_n-\bm \theta)^{\otimes 3}]\right\rangle+\frac{1}{24}\left\langle \nabla^4 R_{\bm\theta}(\bm\theta),\E[(\hbth_n-\bm \theta)^{\otimes 4}]\right\rangle\\
    & + \frac{1}{2n^2}\left\langle H_\theta(\bm\theta), \Lambda(2\Sigma_\theta+M)+\Lambda \bm\theta \bm\theta^\top \Lambda^\top\right\rangle+o(n^{-2}).
\end{align*}
By setting $\Lambda =0$, we recover the expected objective of PTO:
\begin{align*}
     \E[R_{\bm \theta}(\hbth_n )]&= R(\bm \theta)+\frac{1}{2}\left\langle \Htheta, \Sigma_n\right\rangle+\frac{1}{6}\left\langle \nabla^3 R_{\bm\theta}(\bm\theta),\E[(\hbth_n-\bm \theta)^{\otimes 3}]\right\rangle\\
     &+\frac{1}{24}\left\langle \nabla^4 R_{\bm\theta}(\bm\theta),\E[(\hbth_n-\bm \theta)^{\otimes 4}]\right\rangle+o(n^{-2}).
\end{align*}
The improvement is the difference of the last two expressions: 
\begin{align*}
    \E[R_{\bm \theta}((\I+ \Lambda/n)\hbth_n )]-\E[R_{\bm \theta}(\hbth_n )] =  \frac{1}{2n^2}\left\langle H_\theta(\bm\theta), \Lambda(2\Sigma_\theta+M)+\Lambda \bm\theta \bm\theta^\top \Lambda^\top\right\rangle+o(n^{-2}).
\end{align*}
To maximize the difference, note that we are facing a quadratic function of $\Lambda$.
Next we use the property that $\Lambda=\diag(\blam)$ and rewrite it in the standard form $\blam^\top A\blam +b^\top\blam$ for some matrix $A$ and vector $b$.
For the quadratic term $\left\langle H_\theta(\bm\theta), \Lambda \bm\theta \bm\theta^\top \Lambda^\top\right\rangle$, it is easy to see that the coefficient of $\lambda_i\lambda_j$ is $\Htheta_{ij}\theta_i\theta_j$.
Therefore, it can be written as
\begin{equation*}
    \left\langle H_\theta(\bm\theta), \Lambda \bm\theta \bm\theta^\top \Lambda^\top\right\rangle = \blam^\top A(\bm \theta) \blam
\end{equation*}
where $A(\bm \theta)_{ij} = \Htheta_{ij}\theta_i\theta_j$.
For the linear term $\left\langle H_\theta(\bm\theta), \Lambda(2\Sigma_\theta+M)\right\rangle$, we notice that 
the coefficient for $\lambda_i$ is 
\begin{equation*}
    \sum_{j=1}^m \Htheta_{ij}(2\Sigma_\theta+M)_{ij}.
\end{equation*}
As a result, we can write the linear term as $b(\bth)^\top\blam$, where 
\begin{equation*}
    b(\bth) = (\Htheta\odot (2\Sigma_\theta+M) )\one_m,
\end{equation*}
where recall that $\odot$ represents the Hadamard (element-wise) product of two matrices and $\one_m$ is an $m$-dimensional vector of all ones.
This proves the form of the improvement gap \eqref{eq:multi-improvement-gap} in the proposition.
To maximize the gap, the first-order condition yields the claim when $A(\bth)$ is negative definite.
This concludes the proof.
\end{proof}

 \end{APPENDICES}

\end{document}